\documentclass[aos,MSNbibl,dvips]{arximspdf}
\usepackage{mathbh}

\doi{10.1214/14-AOS1246} %kopijuoti is PTS
\volume{42}
\issue{5}
\pubyear{2014}
\firstpage{1941}
\lastpage{1969}
\docsubty{FLA}

\makeatletter
\newcommand{\cT}{\mathcal{T}}
\newcommand{\rn}{\sqrt{n}}
\newcommand{\psg}{{\langle}}
\newcommand{\psd}{{\rangle}}
\newcommand{\sil}{\sigma_{l}}
\newcommand{\eqref}[1]{(\ref{#1})}
\newcommand{\rrvert}{\vert}
\newcommand{\llvert}{\vert}
\newtheorem{teo}{Theorem}
\newtheorem{prop}{Proposition}
\newproclaim{condition}{Condition}
\newproclaim{definition}{Definition}
\newtheorem{corollary}{Corollary}
\newtheorem{lem}{Lemma}
\newproclaim{remark}{Remark}
\def\1{\mathbh{1}}
\newcommand{\given}{\vert}
\makeatother

\begin{document}
\begin{frontmatter}

\title{On the Bernstein--von Mises phenomenon for nonparametric Bayes
procedures}
\runtitle{Nonparametric BvM's}

\begin{aug}
\author[a]{\fnms{Isma\"el}~\snm{Castillo}\ead[label=e1]{ismael.castillo@upmc.fr}\thanksref{t1}}
\and
\author[b]{\fnms{Richard}~\snm{Nickl}\corref{}\ead[label=e2]{r.nickl@statslab.cam.ac.uk}}
\runauthor{I. Castillo and R. Nickl}
\affiliation{CNRS and University of Cambridge}
\address[a]{CNRS---Laboratoire Probabilit\'es\\
\quad et Mod\`eles Al\'eatoires\\
%LPMA - UMR 7599\\
Universit\'es Paris VI and VII\\
B\^atiment Sophie Germain\\
75205 Paris Cedex 13\\
France\\
\printead{e1}} %adresu isvedimo komanda gale!
\address[b]{Statistical Laboratory\\
Department of Pure Mathematics\\
\quad and Mathematical Statistics\\
University of Cambridge\\
CB3 0WB Cambridge\\
United Kingdom\\
\printead{e2}}
\end{aug}
\thankstext{t1}{Supported in part by ANR Grants ``Banhdits'' ANR-2010-BLAN-0113-03 and ``Calibration'' ANR-2011-BS01-010-01.}

% HISTORY:
\received{\smonth{10} \syear{2013}}
\revised{\smonth{3} \syear{2014}}

% ABSTRACT
%
\begin{abstract}
We continue the investigation of Bernstein--von Mises theorems for
nonparametric Bayes procedures from
[\textit{Ann. Statist.} \textbf{41} (2013) 1999--2028]. We introduce
multiscale spaces on which nonparametric priors and posteriors are
naturally defined, and prove Bernstein--von Mises theorems for a
variety of priors in the setting of Gaussian nonparametric regression
and in the i.i.d. sampling model. From these results we deduce several
applications where posterior-based inference coincides with efficient
frequentist procedures, including Donsker-- and Kolmogorov--Smirnov
theorems for the random posterior cumulative distribution functions. We
also show that multiscale posterior credible bands for the regression
or density function are optimal frequentist confidence bands.
\end{abstract}

% KEYWORDS
\begin{keyword}[class=AMS]
\kwd[Primary ]{62G20}
\kwd[; secondary ]{62G15}
\kwd{62G08}
\end{keyword}
\begin{keyword}
\kwd{Bayesian inference}
\kwd{posterior asymptotics}
\kwd{multiscale statistics}
\end{keyword}
\end{frontmatter}

%s1 #&#
\section{Introduction}\label{sec1}

The Bernstein--von Mises (BvM) theorem constitutes a powerful and
precise tool to study Bayes procedures from a frequentist point of
view. It gives universal conditions on the prior under which the
posterior distribution has the approximate shape of a normal
distribution. The theorem is well understood in finite-dimensional
models (see \cite{LeCam86} and \cite{aad98}), but involves some
delicate conceptual and mathematical issues in the infinite-dimensional
setting. There exists a Donsker-type BvM theorem for the cumulative
distribution function based on Dirichlet process priors, see Lo \cite
{L83}, and this carries over to a variety of closely related
nonparametric situations, including quantile inference and censoring
models, where Bernstein--von Mises results are available: see \cite
{CO99,CO04,KL04,K06,HP07} and \cite{HW09}. The
proofs of these results rely on a direct analysis of the posterior
distribution, which is explicitly given in these settings (and
typically of Dirichlet form).

When considering general priors that model potentially smoother
nonparametric objects such as densities or regression functions, the
BvM phenomenon appears to be much less well understood. Notably,
Freedman \cite{Free99} has shown that in a basic Gaussian conjugate
$\ell_2$-sequence space setting, the BvM theorem does not hold true in
generality; see also the related recent contributions \cite{J10,leahu11}. In contrast, in the recent paper \cite{CN13},
nonparametric BvM theorems have been proved in a topology that is
weaker than the one of $\ell_2$, and it was shown that such results
can be useful for several nonparametric problems, including the $\ell
_2$-setting, when applied with care. An important consequence is that,
in contrast to the finite-dimensional situation, whether a
nonparametric posterior credible set is a frequentist confidence set or
not depends in a possibly quite subtle way on the \textit{geometry} of
the set.

The results in \cite{CN13} are confined to the most basic
nonparametric model---Gaussian white noise---and strongly rely on
Hilbert space techniques. The main novelties of the present paper are:
(a) extensions of the results in \cite{CN13} to the i.i.d. sampling
model and (b) the derivation of sharp Bernstein--von Mises results in
spaces whose geometry resembles an $\ell_\infty$-type space and whose
norms are strong enough to allow one to deduce some fundamental new
applications to posterior credible bands and Kolmogorov--Smirnov type
results. Our results are based on mathematical tools developed recently
in Bayesian nonparametrics, particularly the papers \cite{ic13,cr13} and also \cite{rr12}. These give sub-Gaussian estimates on fixed
(semiparametric) functionals of posterior distributions over
well-chosen events in the support of the posterior, which in turn can
be used to control the supremum-type norms relevant in our context via
concentration properties of maxima of sub-Gaussian variables.

Let us outline some applications of our results: consider a prior
distribution $\Pi$ on a family $\mathcal F$ of probability densities
$f$, such as a random Dirichlet histogram or a Gaussian series prior on
the log-density. Let $\Pi(\cdot|X_1,\ldots, X_n)$ be the posterior
distribution obtained from observing $X_1,\ldots, X_n \sim^{\mathrm{i.i.d.}} f$.
%and let $f$ be a draw from $\Pi(\cdot|X_1,\ldots, X_n)$ conditional on
%the $X_i$'s.
It is of interest to study the induced posterior distribution on the
cumulative distribution function $F$ of $f$. Making the ``frequentist''
assumption $X_i \sim^{\mathrm{i.i.d.}}P_0$, the stochastic fluctuations of $F$
around the empirical distribution function $F_n(\cdot) = (1/n)\sum_{i=1}^n 1_{[0,\cdot]}(X_i)$ under the posterior distribution will be
shown to be approximately those of a $P_0$-Brownian bridge $G_{P_0}$:
under the law $P_0^\mathbb N$ of $(X_1, X_2,\ldots)$ the
distributional approximation ($n \to\infty$)
%
%
%e1 #&#
\begin{equation}
\label{donbvm} \sqrt n (F - F_n) \given X_1,\ldots,
X_n \approx G_{P_0}
\end{equation}
holds true, in a sense to be made fully precise below (Corollary~\ref
{classd}). This parallels Lo's \cite{L83} results for the Dirichlet
process and can be used to validate Bayesian Kolmogorov--Smirnov tests
and credible bands from a frequentist point of view. Note, however,
that unlike the results in \cite{L83}, our techniques are not at all
based on any conjugate analysis and open the door to the derivation of
Bernstein--von Mises results in general settings of Bayesian
nonparametrics. We also note that (\ref{donbvm}) is comparable to
central limit theorems $\sqrt n (F_n^b-F_n) \to G_{P_0}$ in
$P_{0}^\mathbb N$-probability for bootstrapped empirical measures
$F_n^b$; see the classical paper \cite{GZ90}. This illustrates how BvM
theorems are in some sense the Bayesian versions of bootstrap
consistency results.

Our results also have important applications for inference on the more
difficult functional parameter $f$ itself. For instance, we will show
that certain $1-\alpha$ posterior credible sets for a density or
regression function are also frequentist optimal, asymptotically exact
level $1-\alpha$ confidence bands.

Before we explain these applications in detail it is convenient to shed
some more light on our general setting. The spaces in which we derive
BvM-type results are in principle abstract and dictated by the
applications we have in mind. They are, however, connected to the
frequentist literature on nonparametric multiscale inference, as
developed in the papers \cite{DS01,DK01,DW08,DKM09,SHMD13}, where also many further references can be
found. This connection gives a further motivation
for our general setting as well as heuristics for the inference
procedures we suggest here. Let us thus explain some main ideas behind
the multiscale approach in the simple regression framework of observing
a signal in Gaussian white noise
%
%
%e2 #&#
\begin{equation}
\label{wnmodel} dX^{(n)}(t) = f(t)\,dt + \frac{1}{\sqrt{n}}\,dW(t),\qquad t
\in[0,1], n \in\mathbb N,
\end{equation}
which can also be written $\mathbb{X}^{(n)}=f+\mathbb{W}/\sqrt{n}$,
with $\mathbb{W}$ a standard white noise; see (\ref{WNmodel}) below
for details. The i.i.d. sampling model, which will be treated below,
gives rise to similar intuitions after replacing $\mathbb{X}^{(n)}-f$
by $P_n-P$ where $P_n = (1/n) \sum_{i=1}^n \delta_{X_i}$ is the
empirical measure from a sample from law $P$ with density~$f$. One
introduces a double-indexed family of linear multiscale functionals
\[
f \mapsto2^{l/2}\int_0^1 \psi
\bigl(2^lx-k\bigr)f(x)\,dx \equiv\langle f, \psi _{lk}
\rangle,
\]
where $l$ is a scaling parameter which has $O(2^l)$ associated location
indices $k$. The prototypical example that we will focus on is to take
a Haar wavelet $\psi= \1_{(0,1/2]}-\1_{(1/2,1]}$, or a more general
wavelet function $\psi$ generating a frame or orthonormal basis $\{
\psi_{lk}\}$ of $L^2$. The projection of $\mathbb{X}^{(n)}-f$ onto
the first $\le J$ scales gives rise to random variables
\[
\sqrt n \bigl\langle \mathbb{X}^{(n)}-f, \psi_{lk} \bigr
\rangle= \langle\mathbb{W}, \psi _{lk} \rangle\equiv g_{lk}
\sim N(0,1),\qquad k,l \le J,
\]
and the maximum over all these statistics scaled by $\sqrt{l}$
%
%
%e3 #&#
\begin{equation}
\label{nulld} Z_J \equiv\sqrt n \max_{l \le J, k}
\frac{|\langle\mathbb
{X}^{(n)}-f, \psi_{lk} \rangle|}{\sqrt l} = \max_{l \le J, k}\frac
{|g_{lk}|}{\sqrt l},
\end{equation}
has a canonical distribution under the null hypothesis $H_0=\{f\}$. The
quantity $Z_J$ is often called a \textit{multiscale statistic}, and
the quantiles of its distribution are used to test hypotheses on $f$.
One can also construct confidence sets $C_n$ by simply taking $C_n$ to
consist of all those $f$ that satisfy simultaneously all the linear constraints
\[
\frac{|\langle\mathbb{X}^{(n)}-f, \psi_{lk} \rangle|}{\sqrt l} \le c_n\qquad\forall k,l,
\]
where $c_n$ are suitable constants chosen in dependence of the
distribution of $Z_J$. Intersecting these linear restrictions with
further qualitative information about $f$, such as smoothness or shape
constraints, can be shown to give optimal frequentist confidence sets
(as, e.g., in Propositions~\ref{histocov} and~\ref{histo0} below).

A key challenge in the multiscale approach is of course the analysis
of the distribution of the random variables $Z_J$. One approach is to
re-center $Z_J$ by a quantity of order $\sqrt J$ and to use extreme
value theory to obtain a Gumbel approximation of the distribution of
these random variables. The slow convergence rates (as $J \to\infty$)
of such limit theorems are often not satisfactory; see, for example,
\cite{H79}. Instead we shall introduce certain sequence spaces in
which direct Gaussian asymptotics can be obtained for multiscale
statistics (without re-centring). This allows for faster convergence
rates (by using standard Berry--Esseen bounds for the central limit
theorem). It is also naturally compatible with a Bayesian approach to
multiscale inference: one distributes independent random variables
across the scales $l$ and locations $k$, corresponding to a random
series prior common in Bayesian nonparametrics. The posterior
distribution then allows one effectively to ``bootstrap'' the law of
$Z_J$, and our BvM-results in multiscale spaces will give a full
frequentist justification of this approach.

Let us illustrate the last point in a key example involving a histogram
prior $\Pi_L, L \in\mathbb N$, equal to the law of the random
probability density
%
%
%e4 #&#
\begin{eqnarray}\label{histo0}
&& f \sim\sum_{k=0}^{2^L-1}
h_k \1_{I_k^L},\qquad I_0^L=\bigl[0,
2^{-L}\bigr], I_k^L = \bigl(k2^{-L},
(k+1)2^{-L}\bigr], k\ge1,
\end{eqnarray}
where the $h_k$ are drawn from a $\mathcal D(1,\ldots, 1)$-Dirichlet
distribution on the unit simplex of $\mathbb R^{2^L}$. Let $\Pi(\cdot
|X_1,\ldots, X_n)$ denote the resulting posterior distribution based
on observing $X_1,\ldots, X_n$ i.i.d. from density $f$. For any
sequence $(w_l)$ such that $w_l/\sqrt l \uparrow\infty$ as $l \to
\infty$ and for standard Haar wavelets
\[
\psi_{-10}=\1_{[0,1]}, \qquad\psi_{lk}=
2^{l/2}(\1_{(k/2^{l}, (k+1/2)/2^{l} ]}-\1_{ ((k+1/2)/2^{l}, (k+1)/2^{l} ]}),
\]
with indices $l \in\mathbb N \cup\{-1,0\}, k =0,\ldots, 2^l-1$, define
%
%
%e5 #&#
\begin{eqnarray}
\label{histo} C_n &\equiv& \biggl\{f\dvtx  \max_{k,l \le L}
\frac{|\langle f- P_n,
\psi_{lk} \rangle|}{w_l} \le\frac{R_n}{\sqrt n} \biggr\},
\end{eqnarray}
where $\langle P_n, \psi_{lk} \rangle= n^{-1} \sum_{i=1}^n \psi
_{lk}(X_i)$ are the empirical wavelet coefficients and where
$R_n=R(\alpha, X_1,\ldots, X_n)$ are random constants chosen such that
\[
\Pi(C_n|X_1,\ldots, X_n) = 1-\alpha,\qquad0<
\alpha<1.
\]

Any set $C_n$ satisfying the identity in the last display is a \textit
{posterior credible set} of level $1-\alpha$, or simply a $(1-\alpha
)$-credible set. Note that in this example the posterior distribution,
and hence $R_n$, can be explicitly computed due to conjugacy of the
Dirichlet distribution under multinomial sampling (i.e., counting
observation points in dyadic bins $I_k^L$).

%The particular credible set $C_n$ in (\ref{histo}) arises from
%simultaneous credible intervals for each scale coefficient $\langle
%f-P_n, \psi_{lk} \rangle$ of the posterior dyadic histogram, with
%critical values $R_n w_l/\sqrt{n}$ at scale $l$, and $R_n$ chosen in
%dependence of the posterior quantiles.

%
%pr1 #&#
\begin{prop} \label{histocov}
Consider the random histogram prior $\Pi$ from (\ref{histo0}) where
$L = L_n$ is such that $2^{L_n} \sim (n/\log n)^{1/2(\gamma+1)}$.
Let $C_n$ be as in (\ref
{histo}). Suppose $X_1,\ldots, X_n$ are i.i.d. from law $P_{0}$ with
density $f_0$ satisfying the H\"older condition
\[
\sup_{x,y \in[0,1], x \neq y}\frac{|f_0(x)-f_0(y)|}{|x-y|^\gamma} <\infty, \qquad1/2<\gamma\le1.
\]
Then we have as $n \to\infty$,
%
%
%e6 #&#
\begin{equation}
\label{histoc} P^\mathbb N_{0}(f_0 \in
C_n) \to1- \alpha.
\end{equation}
Moreover, if $u_n = w_{L_n}/\sqrt{L_n}$, then (\ref{histoc}) remains
true with $C_n$ replaced by
\[
\bar C_n = C_n \cap \bigl\{f\dvtx  \bigl|\langle f,
\psi_{lk} \rangle\bigr| \le u_n 2^{-l (\gamma+1/2)}\ \forall k,l
\bigr\},
\]
and the diameter
$|\bar C_n|_\infty= \sup\{\|f-g\|_\infty\dvtx  f,g \in\bar C_n\}$, satisfies
%
%
%e7 #&#
\begin{equation}
\label{histdiam} |\bar C_n|_\infty=O_{P_{0}^\mathbb N} \biggl(
\biggl(\frac{\log n}{n} \biggr)^{\gamma/(2\gamma+1)} u_n \biggr).
\end{equation}
\end{prop}

We conclude that the $(1-\alpha)$-credible set $C_n$ is an exact
asymptotic frequentist $(1-\alpha)$-confidence set. Following the
multiscale approach, the same is true for~$\bar C_n$ obtained from
intersecting $C_n$ with a $\gamma$-H\"older constraint (expressed
through the decay of the Haar wavelet coefficients). The $L^\infty
$-diameter of $\bar C_n$ shrinks at the optimal rate if the true
density $f_0$ is also $\gamma$-H\"older (noting $u_n \to\infty$ as
slowly as desired). For the proof see Section~\ref{confsec}.

A summary of this article is as follows: in the next section we
introduce the multiscale framework and the statistical sampling models
and show how to construct efficient frequentist estimators in them. In
Section~\ref{BAY} we introduce the Bayesian approach, formulate a
general notion of a nonparametric Bernstein--von Mises phenomenon in
multiscale spaces and prove that the phenomenon occurs for a variety
of relevant nonparametric prior distributions, including Gaussian
series priors and random histograms. In Section~\ref{app} we discuss
statistical applications to Donsker--Kolmogorov--Smirnov theorems and
credible bands. Section~\ref{proof} contains the proofs.

%s2 #&#
\section{The general framework}

We use the usual notation for $L^p=L^p([0,1])$-spaces of integrable
functions, and we denote by $\ell_p$ the usual sequence spaces. The
usual supremum norm is denoted by $\|\cdot\|_\infty$. Throughout we
consider an $S$-regular, $S \ge0$, wavelet basis
%
%
%e8 #&#
\begin{equation}
\label{wavonb} \bigl\{\psi_{lk}\dvtx  l \ge J_0-1, k=0,\ldots, 2^{l}-1\bigr\}, \qquad J_0 \in \mathbb N \cup\{0
\},
\end{equation}
of $L^2([0,1])$ (by convention we denote the usual ``scaling function''
$\varphi$ as the first wavelet $\psi_{(J_0-1)0}$). We restrict to
Haar wavelets ($S=J_0=0$), periodised wavelet bases ($J_0=0, S>0$) or
boundary corrected wavelet bases ($S>0$, $J_0=J_0(S)$ large enough, see
\cite{CDV93}). Functions $f \in L^2$ generate double-indexed sequences
$\{\langle f, \psi_{lk}\rangle= \int_0^1 f \psi_{lk}\}$, and
conversely any sequence $(x_{lk})$ generates wavelet series of
(possibly generalised) functions $\sum_{k,l} x_{lk} \psi_{lk}$ on $[0,1]$.

We define H\"{o}lder-type spaces $C^s$ of continuous functions on $[0,1]$:
%
%
%e9 #&#
\begin{equation}
\label{hold} C^s\bigl([0,1]\bigr) = \Bigl\{f \in C\bigl([0,1]
\bigr)\dvtx  \|f\|_{s, \infty}:= \sup_{l, k} 2^{l(s+1/2)} \bigl|
\langle\psi_{lk}, f \rangle\bigr| < \infty \Bigr\}.
\end{equation}
When the wavelets are regular enough, this norm characterises the scale
of H\"{o}lder (--Zygmund when $s \in\mathbb N$) spaces. Otherwise we
work with the spaces defined through decay of the multiscale
coefficients, which still contain the classical \mbox{$s$-}H\"older spaces by
standard results in wavelet theory.

Convergence in distribution of random variables $X_n \to^d X$ in a
metric space $(S,d)$ can be metrised by metrising weak convergence of
the induced laws $\mathcal L(X_n)$ to $\mathcal L(X)$ on $S$. For
convenience we work with the bounded-Lipschitz metric $\beta_S$: let
$\mu, \nu$ be probability measures on $(S,d)$, and define
%
%
%e10 #&#
\begin{eqnarray}\label{blmet}
\beta_S(\mu, \nu) &\equiv&\sup_{F\dvtx  \|F\|_{BL}\le1}
\biggl\llvert \int_S F(x) \bigl(d\mu(x)-d \nu(x)\bigr)
\biggr\rrvert,
\nonumber\\[-8pt]\\[-8pt]\nonumber
\|F\|_{BL}&=& \sup_{x \in S}\bigl|F(x)\bigr| + \sup
_{x\neq y, x, y \in S} \frac
{|F(x)-F(y)|}{d(x,y)}.
\end{eqnarray}

%s2.1 #&#
\subsection{Multiscale spaces}

For monotone increasing weighting sequences $w=(w_l\dvtx  l \ge J_0-1), w_l
\ge1$, we define
multiscale sequence spaces
%
%
%e11 #&#
\begin{equation}
\mathcal M \equiv\mathcal M(w) \equiv \biggl\{x=\{x_{lk}\}\dvtx  \|x\|
_{\mathcal M(w)} \equiv\sup_{l}\frac{\max_{k}|x_{lk}|}{w_l} <\infty
\biggr\}.
\end{equation}
The space $\mathcal M(w)$ is a nonseparable Banach space (it is
isomorphic to $\ell_\infty$). The (weighted) sequences in $\mathcal
M(w)$ that vanish at infinity form a separable closed subspace for the
same norm
%
%
%e12 #&#
\begin{equation}
\label{em0} \mathcal M_0=\mathcal M_0(w)= \biggl
\{x \in\mathcal M(w)\dvtx  \lim_{l \to
\infty} \max_k
\frac{|x_{lk}|}{w_l}=0 \biggr\}.
\end{equation}

We notice that $w_l \ge1$ implies $\|x\|_{\mathcal M} \le\|x\|_{\ell
_2}$ so that $\mathcal M$ always contains $\ell_2$. For suitable
divergent weighting sequences $(w_l)$, these spaces contain objects
that are much less regular than $\ell_2$-sequences, such as a Gaussian
white noise $dW$. The action of $dW$ on $\{\psi_{lk}\}$ generates an
i.i.d. sequence $g_{lk}$ of standard $N(0,1)$'s, hence whether $dW$
defines a Gaussian Borel random variable $\mathbb W$ in $\mathcal M_0$
or not depends entirely on the weighting function $w$.

%
%de1 #&#
\begin{definition}\label{admi}
Call a sequence $(w_l)$ admissible if $w_l/\sqrt l \uparrow\infty$ as
$l \to\infty$.
\end{definition}
%
%
%pr2 #&#
\begin{prop} \label{CKR} Let $\mathbb W = (\int\psi_{lk} \,dW\dvtx  l,k)=
(g_{lk}), g_{lk} \sim N(0,1)$, be a Gaussian white noise. For $\omega
=(\omega_l)=\sqrt l$ we have
$E\|\mathbb W\|_{\mathcal M(\omega)} < \infty$. If $w=(w_l)$ is
admissible, then $\mathbb W$ defines a tight Gaussian Borel probability
measure in the space $\mathcal M_0(w)$.
\end{prop}
\begin{pf}
Since there are $2^l$ i.i.d. standard Gaussians $g_{lk}=\langle\psi
_{lk}, dW \rangle$ at the $l$th level, we have from a standard bound
$E \max_{k}|g_{lk}| \leq C \sqrt l$ for some universal constant $C$.
The Borell--Sudakov--Tsirelson inequality (e.g., \cite{L01}) applied
to the maximum at the $l$th level gives, for any $M$ large enough,% and
%a universal constant $C>0$,
%
\begin{eqnarray*}
&& \Pr \Bigl(\sup_l l^{-1/2} \max
_k |g_{lk}| > M \Bigr)
\\
&&\qquad \le\sum_l \Pr \Bigl( \max_k|g_{lk}|
- E \max_k |g_{lk}| > \sqrt l M - E \max
_k |g_{lk}| \Bigr)
\\
&&\qquad \leq2 \sum_l \exp \bigl\{- c(M-C)^2
l \bigr\}.
\end{eqnarray*}
Now using $E[X] \le K + \int_K^\infty\Pr[X\ge t]\,dt$ for any
real-valued random variable $X$ and any $K\ge0$, one obtains that $\|
\mathbb W\|_{\mathcal M(\omega)}$ has finite expectation.

It now also follows immediately from the definition of the space
$\mathcal M_0(w)$ that for any sequence $w_l/\sqrt l \uparrow\infty$,
we have $\mathbb W \in\mathcal M_0$ almost surely. Since the latter is
a separable complete metric space, $\mathbb W$ is a tight Gaussian
Borel random variable in it (e.g., page 374 in \cite{B98}).
\end{pf}

%
%re1 #&#
\begin{remark}[(Admissible sequences $w$)] \label{CKRop}
Assuming admissibility of $w$ is necessary if one wants to
show that $\mathbb W$ is tight in $\mathcal M(w)$. Since weak
convergence of probability measures on a complete metric space implies
tightness of the limit distribution, it is in particular impossible, as
will be relevant below, to converge weakly towards $\mathbb W$ in
$\mathcal M(w)$ without assuming admissibility of $w$. To prove that
admissibility is necessary, suppose on the contrary that $\mathbb W$
were tight in $\mathcal M(\omega)$ for some sequence $\omega_l \sim
\sqrt l$, hence defining a Radon Gaussian measure in that space. Then
by Theorem 3.6.1 in \cite{B98} the topological support of $\mathbb W$
equals the completion of the RKHS $\ell_2$ in the norm of the ambient
Banach space $\mathcal M(\omega)$, which is $\mathcal M_0(\omega)$. Since
\[
\lim_{J \to\infty} \frac{\max_k |g_{Jk}|}{\sqrt J} = \sqrt{2 \log 2} \neq0
\]
almost surely
%book- and next use Borel-Cantelli ?}
we have $\mathbb W \notin\mathcal M_0(\omega)$, a contradiction, so
$\mathbb W$ cannot be tight. The cylindrically-defined law of $\mathbb
W$ is in fact a ``degenerate'' Gaussian measure in $\mathcal M(\omega
)$ that does (assuming the continuum hypothesis) not admit an extension
to a Borel measure on $\mathcal M(\omega)$; see Definition 3.6.2 and
Proposition 3.11.5 in~\cite{B98}. It has further unusual properties:
$\mathbb W$ has a ``hole.'' That is, for some $c>0$, $\|\mathbb W\|
_{\mathcal M(\omega)} \in[c, \infty)$ almost surely (see \cite
{CKR93}), and depending on finer properties of the sequence $\omega$,
the distribution of $\|\mathbb W\|_\mathcal M$ may not be absolutely
continuous, and its absolutely continuous part may have infinitely many
modes; see \cite{HJDS79}.
\end{remark}

% \begin{remark} For $w_l=1$ the space $\mathcal M$ in fact coincides
%with the wavelet coefficient sequence spaces arising from elements of
%negative order Besov spaces $ B^{-1/2}_{\infty\infty}([0,1])$, which
%partly explains the relationship of this approach to the one in
%has a weaker norm than the spaces $\ell^\infty(\mathcal G)$ for $
%have in mind the multiscale norm is however strong enough.

%s2.2 #&#
\subsection{Nonparametric statistical models}

%s2.2.1 #&#
\subsubsection{Nonparametric regression}
For $f \in L^2$ consider observing a trajectory in the white noise
model (\ref{wnmodel}) which is a natural surrogate for a fixed design\vadjust{\goodbreak}
nonparametric regression model with Gaussian errors. By Proposition
\ref{CKR} and since any $f \in L^2$ has wavelet coefficients $\{
f_{lk}\} \in\ell_2\subset\mathcal M_0(w)$, equation (\ref{wnmodel})
makes rigorous sense as the tight Gaussian shift experiment
%
%
%e13 #&#
\begin{equation}
\label{WNmodel} \mathbb X^{(n)}= f + \frac{1}{\sqrt n} \mathbb W,\qquad
n \in\mathbb N,
\end{equation}
in $\mathcal M_0(w)$ for any admissible $(w_l)$. We denote the law
$\mathcal L(\mathbb{X}^{(n)})$
by $P^n_{f}$. Then
%
%
%e14 #&#
\begin{equation}
\label{gausdon} \sqrt n \bigl(\mathbb X^{(n)}-f\bigr) = \mathbb W\qquad
\mbox{in } \mathcal M_0,
\end{equation}
and one deduces that $\mathbb{X}^{(n)}$ is an efficient estimator of
$f$ in $\mathcal M_0$.

%s2.2.2 #&#
\subsubsection{The i.i.d. sampling setting}

Consider next the situation where we observe $X_1,\ldots, X_n$
i.i.d. from law $P$ with density $f$ on $[0,1]$. Then a natural
estimate of $\langle f, \psi_{lk} \rangle$ is given by $P_n\psi
_{lk}\equiv\langle P_n,\psi_{lk}\rangle= \frac{1}{n} \sum_{i=1}^n
\psi_{lk}(X_i)$. By the central limit theorem, for $k,l$ fixed and as
$n \to\infty$, the random variable $\sqrt{n}(P_n - P)(\psi_{lk})$
converges in distribution to
%
%
%e15 #&#
\begin{equation}
\label{GPwhite} \mathbb G_{P}(\psi_{lk}) \sim N\bigl(0,
\operatorname{Var}_P\bigl(\psi_{lk}(X_1)
\bigr)\bigr).
\end{equation}
In analogy to the white noise process $\mathbb{W}$, the process
$\mathbb G_P$ arising from (\ref{GPwhite}) can be rigorously defined
as the Gaussian process indexed by the Hilbert space
\[
L^2(P)\equiv \biggl\{f\dvtx [0,1]\to\mathbb R\dvtx  \int_0^1
f^2\,dP<\infty \biggr\}
\]
with covariance function
%isonormal process on $L^2_0(P):=\{g\in L^2(P), \int gdP=0\}$, equipped
%with the inner-product
$\mathbb{E} [\mathbb G_{P}(g) \mathbb G_{P}(h) ]=\int_0^1
(g-Pg)(h-Ph)\,dP$. We call $\mathbb G_P$ the \emph{$P$-white bridge
process}. An analogue of Proposition~\ref{CKR}, and of the remark
after it, holds true for $\mathbb G_P$ whenever $P$ has a bounded density.

%If $P$ is the uniform distribution then the law of $\mathbb G_{P}$ in $
%analogue of Proposition~\ref{CKR}, and of the remark after it, holds
%true whenever $P$ has a bounded density.

%
%pr3 #&#
\begin{prop} \label{CKRP} Proposition~\ref{CKR} holds true for the
$P$-white bridge $\mathbb G_{P}$ replacing $\mathbb W$ whenever $P$ has
a bounded density on $[0,1]$.
\end{prop}
\begin{pf}
The proof is exactly the same, using the standard bounds
\[
\operatorname{Var}\bigl(\mathbb G_{P}(\psi_{lk})\bigr)
\le\|f\|_\infty,\qquad E \max_{k}\bigl|\mathbb
G_{P}(\psi_{lk})\bigr| \le C \|f\|^{1/2}_\infty
\sqrt l,
\]
where $f$ denotes the density of $P$.
\end{pf}

Any $P$ with bounded density $f$ has coefficients $\langle f, \psi
_{lk} \rangle\in\ell_2 \subset\mathcal M_0(w)$. We would like to
formulate a statement such as
\[
\sqrt n (P_n -P) \to^d \mathbb G_P\qquad
\mbox{in } \mathcal M_0,
\]
as $n \to\infty$, paralleling (\ref{gausdon}) in the Gaussian white
noise setting. The fluctuations of $\sqrt n (P_n-P)(\psi_{lk})/\sqrt
{l}$ along $k$ are stochastically bounded for $l$ such that $2^l \le
n$, but are unbounded for high frequencies. Thus the empirical process
$\sqrt n (P_n-P)$ will not define an element of $\mathcal M_0$ for
every admissible sequence $w$. In our nonparametric setting we can
restrict to frequencies at levels $l, 2^l \le n$ and introduce an
appropriate ``projection'' $P_n(j)$ of the empirical measure $P_n$ onto
$V_j$ via
%
%
%e16 #&#
\begin{equation}
\label{def-pn} \bigl\langle P_n(j), \psi_{lk}\bigr
\rangle= \cases{ {\langle}P_n, \psi_{lk} {\rangle}, &\quad
if $l\le j$,
\cr
0, &\quad if $l> j$,}
\end{equation}
which defines a tight random variable in $\mathcal M_0$. The following
theorem shows that $P_n(j)$ estimates $P$ efficiently in $\mathcal M_0$
if $j$ is chosen appropriately. Note that the natural choice $j=L_n$
such that
\[
2^{L_n} \sim N^{1/(2\gamma+1)},
\]
where $N=n$ (if $\gamma>0$) or $N=n/\log n$ (if $\gamma\ge0$), is possible.

%
%th1 #&#
\begin{teo} \label{CLT}
Let $w=(w_l)$ be admissible. Suppose $P$ has density $f$ in $C^\gamma
([0,1])$ for some $\gamma\ge0$. Let $j_n$ be such that
\[
\sqrt n 2^{-j_n(\gamma+1/2)}w^{-1}_{j_n} = o(1),\qquad
\frac{2^{j_n}
j_n}{n} =O(1).
\]
Then we have, as $n \to\infty$,
\[
\sqrt n \bigl(P_n(j_n) - P\bigr) \to^d
\mathbb G_P\qquad\mbox{in } \mathcal M_0(w).
\]
\end{teo}

%s3 #&#
\section{The nonparametric Bayes approach} \label{BAY}

In both regression or density estimation one constructs a prior
probability distribution from which the function $f$ is drawn, and
given the observations\vadjust{\goodbreak} $X=X^{(n)}$, equal to\vspace*{1pt} either $\mathbb
{X}^{(n)}\sim P^n_f$ or $X_1,\ldots, X_n$ i.i.d. from density $f$, one
computes the posterior distribution $\Pi(\cdot|X)$ of $f$. Under
appropriate conditions the wavelet coefficient sequence associated to a
posterior draw $f \sim\Pi(\cdot|X)$ will give rise to a random
variable in $\mathcal M_0$. If $T_n=T_n(X)$ is an efficient estimator
of $f$ in $\mathcal M_0$, such as $\mathbb X^{(n)}$ or $P_n(j)$ from
the previous subsections, then one can ask, following \cite{CN13}, for
a Bernstein--von Mises type result: assuming $X \sim P_{f_0}$ for some
fixed $f_0$, do we have
%as $n \to\infty$,
%
%e17 #&#
\begin{equation}
{\mathcal{L}}\bigl(\sqrt n (f- T_n) \given X\bigr) \to\mathcal L(
\mathbb G)\qquad\mbox{weakly in } \mathcal M_0(w)\mbox{ as }n
\to \infty,
\end{equation}
with $P_{f_0}$-probability close to one? Here, depending on the
sampling model considered, $\mathbb G$ equals either $\mathbb W$ or
$\mathbb G_{P_0}, dP_0(x) = f_0(x)\,dx$ and $P_{f_0}$ stands, in slight
abuse of notation, for the law $P_{f_0}^n$ of $\mathbb{X}^{(n)}$ or
the law $P_{0}^\mathbb N$ of $(X_1, X_2,\ldots)$.

To make such a statement rigorous we will metrise weak convergence of
laws in $\mathcal M_0(w)$ via $\beta_{\mathcal M_0(w)}$ from (\ref
{blmet}), and view the prior $\Pi$ on the functional parameter $f\in
L^2$ as a prior on sequence space $\ell_2$ under the wavelet isometry
$L^2 \cong\ell_2$ [arising from an arbitrary but fixed wavelet basis
(\ref{wavonb})].

%
%de2 #&#
\begin{definition}\label{weakbvmp}
Let $w$ be admissible, let $\Pi$ be a prior and $\Pi(\cdot|X)$ the
corresponding posterior distribution on $\ell_2 \subset\mathcal
M_0=\mathcal M_0(w)$, obtained from observations $X$ in the white noise
or i.i.d. sampling model. Let $\tilde\Pi_n$ be the image measure
%%obtained from
of $\Pi(\cdot|X)$ under the mapping
\[
\tau\dvtx  f \mapsto\sqrt n (f-T_n),
\]
where $T_n=T_n(X)$ is an estimator of $f$ in $\mathcal M_0$. Then we
say that $\Pi$ satisfies the weak Bernstein--von Mises phenomenon in
$\mathcal M_0$ with centring $T_n$ if, for $X \sim P_{f_0}$ and fixed
$f_0$, as $n \to\infty$,
\[
\beta_{\mathcal M_0}(\tilde\Pi_n, \mathcal N) \to^{P_{f_0}} 0,
\]
where $\mathcal N$ is the law in $\mathcal M_0$ of $ \mathbb W$ or of
$ \mathbb G_{P_{0}}, f_0 \in L^\infty$, respectively.
\end{definition}

%
%re2 #&#
\begin{remark}
If convergence of moments (Bochner-integrals) $E[\tilde\Pi_n|X] \to
^{P_{f_0}} E\mathcal N=0$ occurs in the above limit, then we deduce
%
%
%e18 #&#
\begin{equation}
\label{mom} \|\bar f_n - T_n\|_{\mathcal M_0} =
o_{P_{f_0}}(1/\sqrt n),
\end{equation}
where\vspace*{1pt} $\bar f_n=E(f|X)$ is the posterior mean. If $T_n$ is an efficient
estimator of \mbox{$f \in\mathcal M_0$}, then (\ref{mom}) implies that $\bar
f_n$ is so too.
\end{remark}

In \cite{CN13}, Bernstein--von Mises theorems are proved in certain
negative Sobolev spaces $H(\delta), \delta>1/2$, and various
applications of such results are presented. A~multiscale BvM result in
$\mathcal M_0$ for a prior $\{f_{lk}\}$ implies a weak BvM for the
prior $\sum_{k,l} f_{lk} \psi_{lk}$ in $H(\delta)$, as the following
result shows. In particular all the applications from \cite{CN13}
carry over to the present setting.

%
%pr4 #&#
\begin{prop}
Suppose the weak Bernstein--von Mises phenomenon holds true in
$\mathcal M_0(w)$ with $(w_l)$ such that $\sum_{l} w^2_l l^{-2\delta
}<\infty$ for some $\delta>0$. Then the weak Bernstein--von Mises
phenomenon holds in $H(\delta)$.
\end{prop}
\begin{pf}
The norm of $H(\delta)$ is given by (see \cite{CN13}, Section~1.2),
\begin{eqnarray*}
\|f\|^2_{H(\delta)} &=& \sum_l
2^{-l}l^{-2\delta} \sum_{k} \bigl|\langle
f, \psi_{lk} \rangle\bigr|^2
\\
& \le&\sup_{l} w_l^{-2} \max
_k\bigl|\langle f, \psi_{lk} \rangle\bigr|^2
\sum_{l} w^2_l
l^{-2\delta}
\\
&\le& C\|f\|_{\mathcal M_0(w)}^2,
\end{eqnarray*}
so that the result follows from the continuous mapping theorem.
\end{pf}

While the above notions of the BvM phenomenon will be shown below to be
useful and feasible in nonparametric settings, there are other ways to
formulate BvM-type statements. For instance, one may investigate how the
classical BvM theorem in finite-dimensions extends to parameter spaces
of dimension that increases with $n$; see, \cite{G99,B11,K13} for
results in this direction.

%We now prove that the weak BvM holds true in $\cM_0$ in a variety of
%concrete examples.

Throughout the rest of this section $ {\mathcal{M}}_0= {\mathcal
{M}}_0(w)$ is the space defined in (\ref{em0}), with $w$ an admissible
sequence as in Definition~\ref{admi}.

%s3.1 #&#
\subsection{Bernstein--von Mises theorems in $\mathcal M_0(w)$: Gaussian regression case}

In the white noise model (\ref{WNmodel}) natural priors for $f$ are
obtained from distributing random coefficients on the $\psi_{lk}$'s.

%
%co1 #&#
\begin{condition} \label{tolk}
Consider product priors $\Pi$ arising from random functions
\[
f(x) = \sum_{l} \sigma_l \sum
_{k } \phi_{lk} \psi_{lk}(x), \qquad x
\in[0,1],
\]
where the $\phi_{lk}$ are i.i.d. from probability density $\varphi\dvtx
\mathbb R \to[0, \infty)$ satisfying
{\renewcommand{\theequation}{E}
%e19 #&#
\begin{equation}\label{eqE}
\exists a,C >0\qquad\forall x\in\mathbb{R}, \qquad\varphi(x)\le
Ce^{-a x^2},
\end{equation}}%
and where $\sigma_l = 2^{-l(\alpha+ (1/2))}, \alpha>0$,
ensuring in particular that $f \in L^2$ almost surely.

For $\mathbb X^{(n)} \sim P_{f_0}^n$ and $f_0$ with wavelet
coefficients $\{\langle f_0, \psi_{lk} \rangle\} \in\ell_2$, we
assume moreover that there exists a finite constant $M>0$
such that
{\renewcommand{\theequation}{P1}
%e20 #&#
\begin{equation}\label{eqP1}
\sup_{l,k} \frac{|\langle f_0, \psi_{lk} \rangle
|}{\sigma_{l}} \le M,
\end{equation}}%
and that there exists $\tau>M, c_\varphi>0$ such that
{\renewcommand{\theequation}{P2}
%e21 #&#
\begin{equation}\label{eqP2}
\mbox{on } (-\tau,\tau)\mbox{ the density } \varphi
\mbox{ is continuous and satisfies } \varphi\ge c_\varphi.
\end{equation}}%
\end{condition}\setcounter{equation}{18}
If $f_0 \in C^\beta, \beta> 0$, then (\ref{eqP1}) is satisfied as soon
as $\alpha\le\beta$ (so any prior that matches the regularity of
$f_0$, or that ``undersmooths,'' can be used).

%
%re3 #&#
\begin{remark}
Condition~\ref{tolk} allows for a sub-Gaussian density $\varphi$.
Strictly subexponential tails could be allowed too if the weighting
sequence $w$ satisfies an additional constraint: Theorem~\ref{wnbvm}
below holds true for exponential-power densities $\varphi(x) \asymp
e^{-|x|^p}$ and $1<p<2$, provided $w_l/l^{1/p} \uparrow\infty$.
\end{remark}

Any prior satisfying Condition~\ref{tolk} defines a Borel probability
measure on $L^2$ (using separability of the latter space), and the
resulting posterior distribution also defines an element of $L^2\cong
\ell_2\subset\mathcal M_0$.
%Since $L^2 = \ell_2$ is continuously injected into the multiscale
%space $\mathcal M_0$ we can view the posterior random measure as an
%element
%$\Pi(\cdot|\mathbb X^{(n)})$ in $\mathcal M_0$.

%
%th2 #&#
\begin{teo} \label{wnbvm}
Suppose $\Pi$ satisfies Condition~\ref{tolk}, and let $\Pi(\cdot
|\mathbb X^{(n)})$ be the posterior distribution in $\mathcal M_0$
arising from observing (\ref{WNmodel}) for some fixed $f_0 \in C^\beta
$, $\beta>0$. Then $\Pi$ satisfies the weak Bernstein--von Mises
theorem in the sense of Definition~\ref{weakbvmp} in the space
$\mathcal M_0= \mathcal M_0(w)$ for any admissible $w$, with $\mathcal
N$ equal to the law of $\mathbb W$, and with centring $T_n$ equal to
$\mathbb X^{(n)}$ or equal to the posterior mean $E(f|\mathbb X^{(n)})$.
\end{teo}

%s3.2 #&#
\subsection{Bernstein--von Mises theorems in $\mathcal M_0(w)$: Sampling model case}

Let\vspace*{1pt} us now turn to the situation where one observes a sample $X_i \sim
^{\mathrm{i.i.d.}} P$,
\[
(X_1,\ldots, X_n) \equiv X^{(n)},
\]
from law $P$ with \textit{bounded} probability density $f$ on $[0,1]$.
We define multiscale priors $\Pi$ on some space $\mathcal F$ of
probability density functions $f$ giving rise to absolutely continuous
probability measures. Let
%More precisely, for a given $\rho>0$ and $D<\infty$, first set
%
\[
{\mathcal{F}}:= \bigcup_{0<\rho\le D<\infty} {\mathcal {F}}(
\rho,D):= \bigcup_{0<\rho\le D<\infty} \biggl\{f\dvtx [0,1]\to[\rho, D],
\int_0^1f=1 \biggr\}.
\]
%
%Let $\cF$ be the union of all sets $\cF(\rho,D)$, $0<\rho<D<\infty$.
In the following we assume that the ``true'' density $f_0$ belongs to
$ {\mathcal{F}}_0:= {\mathcal{F}}(\rho_0,D_0)$, for some
$0<\rho_0\le D_0<\infty$.

We consider various classes of priors on densities and two possible
values for a cut-off parameter $L_n$.
For $\alpha>0$, let $j_n=j_n(\alpha)$ and $l_n=l_n(\alpha)$ be the
\textit{largest} integers such that
%
%
%e19 #&#
\begin{equation}
\label{def-both} 2^{j_n} \le n^{1/(2\alpha+1)},\qquad2^{l_n} \le
\biggl(\frac
{n}{\log n} \biggr)^{1/(2\alpha+1)},
\end{equation}
and set, in slight abuse of notation, either
%
%
%e20 #&#
\begin{equation}
\label{def-ln} L_n = j_n\qquad (\forall n\ge1) \quad\mbox{or}
\quad L_n = l_n\qquad (\forall n\ge1).
\end{equation}

(\textbf{S})\label{eqS} \textit{Priors on log-densities}. %We consider two main families
%of priors within this class.
%Some applications will be shown to hold for the three families, while
%for some others we will restrict to the first two families.
%{\bf(S1) - \bf(S2)} {\it Truncated series with heavy tails/with
%Gaussian priors}.
Given a multiscale wavelet basis $\{\psi_{lk}\}$ from (\ref
{wavonb}), consider the prior $\Pi$
induced by, for any $x\in[0,1]$ and $L_n$ as in (\ref{def-ln}),
%
%
%e21 #&#
%e22 #&#
\begin{eqnarray}
T(x) & =& \sum_{l\le L_n} \sum
_{k=0}^{2^l -1} \sigma_{l}
\alpha_{lk} \psi_{lk}(x), \label{def-T}
\\
f(x) & =& \exp \bigl\{ T(x) - c(T) \bigr\},\qquad c(T)=\log\int
_0^1 e^{T(x)} \,dx, \label{def-prior-log}
\end{eqnarray}
where $\alpha_{lk}$ are i.i.d. random variables of continuous
probability density $\varphi\dvtx  \mathbb R \to[0,\infty)$. We consider
the choices
{\renewcommand{\theequation}{\textbf{S\arabic{equation}}}\setcounter{equation}{0}
%e1 #&#
%e2 #&#
\begin{eqnarray}
\label{eqS1}
\varphi(x) & =& \varphi_H(x),
\\
\label{eqS2} \varphi(x) &=& \varphi_G(x)=e^{-x^2/2}/\sqrt{2\pi},
\end{eqnarray}}%
where $\varphi_H$ is any density such that $\log\varphi_H$ is
Lipschitz on $\mathbb{R}$. We call this the log-Lipschitz case. For
instance, the $\alpha_{lk}$'s can be Laplace-distributed or have
heavier tails. To simplify some proofs we restrict to a specific form
of density:
for a given $0\le\tau<1$ and
$x\in\mathbb{R}$, and $c_\tau$ a normalising constant, suppose
$\varphi_H$ takes the form
\setcounter{equation}{22}%e23 #&#
\begin{equation}
\varphi_{H,\tau}(x) = c_\tau\exp\bigl\{-\bigl(1+|x|\bigr)^{1-\tau}
\bigr\}. \label{phi-ht}
\end{equation}
Suppose the prior parameters $\sigma_{l}$ satisfy, for $\alpha>1/2$
and $0<r < \alpha-1/4$,
%
%
%e24 #&#
\begin{eqnarray}\label{cond-sil}
\sigma_{l}&=&2^{-l(\alpha+1/2)}\qquad\mbox{(log-Lipschitz-case)},
\nonumber\\[-8pt]\\[-8pt]
\sigma_{l}&=& 2^{-l(r+(1/2))}\qquad\mbox{(Gaussian-case)}.\nonumber
\end{eqnarray}

(\textbf{H})\label{eqH} \textit{Random histograms density priors}.
Associated to the regular dyadic partition of $[0,1]$ at level $L \in
\mathbb{N} \cup\{0\}$, given by
$I_0^L=[0,2^{-L}]$ and $I_k^L=(k2^{-L},(k+1)2^{-L}]$ for $k=1,\ldots,
2^L-1$, is a natural notion of histogram
\[
{\mathcal{H}}_L = \Biggl\{h\in L^\infty[0,1],  h(x) = \sum
_{k=0}^{2^L-1} h_k
\1_{I_k^L}(x),  h_k \in\mathbb{R}, k=0,\ldots,2^L-1 \Biggr\}.
\]
Let
$\mathcal S_L = \{ \omega\in[0,1]^{2^L}; \sum_{k=0}^{2^L-1} \omega
_k =1\}$
be the unit simplex in $\mathbb R^{2^L}$. Further denote
%2^L \sum_{k=0}^{2^L-1} \omega_k \1_{\ikl}(x),\qquad(\omega_0,\ldots,
$ {\mathcal{H}}_L^1$ the subset of $ {\mathcal{H}}_L$
consisting of histograms which are densities on $[0,1]$ with~$L$
equally spaced dyadic knots. Let $ {\mathcal{H}}^1$ be the set of
all histograms which are densities on $[0,1]$.
%The set $\cH_L$ is a closed subspace of $L^2[0,1]$. For any function
%$h$ in $L^2[0,1]$, consider its projection $h_{[L]}$ in the
%$L^2$-sense on $\cH_L$.
%It holds
% k\sum_{j=1}^k \bar{h}_{j,k} \1_{I_j}. Useful ?
%Useful elementary properties on histograms are gathered in Lemma
%[lem-hist] below.

%For the next Theorem, the following notation is useful. Denote by

A simple way to specify a prior $\Pi$ on $ {\mathcal{H}}_L^1$ is
to set $L=L_n$ deterministic and to fix a distribution
for $\omega_L:=(\omega_0,\ldots,\omega_{2^L-1})$.
Set $L=L_n$ as defined in (\ref{def-ln}).
Choose some fixed constants $a, c_1,c_2>0$, and let
%and set, for $L_n=\lfloor\log_2 J_n \rfloor$ and $J_n= \lfloor(n/
%
%e25 #&#
\begin{equation}
\label{prk} L = L_n,\qquad\omega_L \sim{\mathcal{D}}(
\alpha_0,\ldots, \alpha_{2^L-1}), \qquad c_1
2^{-La}\le\alpha_k \le c_2,
\end{equation}
for any admissible index $k$, where $ {\mathcal{D}}$ denotes the
Dirichlet distribution on $\mathcal S_L$. Unlike those suggested by the
notation, the coefficients $\alpha$ of the Dirichlet distribution are
allowed to depend on $L_n$, so that $\alpha_{k}=\alpha_{k,L_n}$.

The priors {(S)}, {(H)} above are ``multiscale'' priors where high
frequencies are ignored, corresponding to truncated series priors
considered frequently in the nonparametric Bayes literature. The
resulting posterior distributions $\Pi(\cdot|X^{(n)})$ attain minimax
optimal contraction rates up to logarithmic terms in Hellinger and
\mbox{$L^2$-}distance \cite{vvvz,rr12,cr13} and $L^\infty
$-distance \cite{ic13}. Clearly other priors are of interest as
well, for instance, priors without or with random high-frequency
cut-off or Dirichlet mixtures of normals etc. While our current proofs
do not cover such situations, one can note that our proof strategy via
simultaneous control of many linear functionals is applicable in such
situations as well. Generalising the scope of our techniques is an
interesting direction of future research.

The projection $P_n(j)$ as in (\ref{def-pn}), with the choice $j=L_n$
from (\ref{def-ln}), defines a tight random variable\vspace*{2pt} in $\mathcal
M_0$. For $z\in\mathcal M_0$, the map $\tau_z\dvtx  f \mapsto\sqrt n (f-
z)$ maps $\mathcal M_0 \to\mathcal M_0$, and we can define the shifted
posterior $\Pi(\cdot\given X^{(n)}) \circ\tau_{P_n(L_n)}^{-1}$. The
following theorem shows that the above priors satisfy a weak BvM
theorem in $\mathcal M_0$ in the sense of Definition~\ref{weakbvmp},
with efficient centring $P_n(L_n)$; cf. Theorem~\ref{CLT}. Denote the
law $\mathcal L(\mathbb G_{P_0})$ of $\mathbb G_{P_0}$ from Proposition
\ref{CKRP} by $\mathcal N$.

%
%th3 #&#
\begin{teo}\label{thmdens1}
Let $\mathcal M_0=\mathcal M_0(w)$ for any admissible $w=(w_l)$. Let
$X^{(n)}=(X_1,\ldots, X_n)$ i.i.d. from law $P_0$ with density $f_0\in
{\mathcal{F}}_0$. Let $\Pi$ be a prior on the set of probability
densities $\mathcal F$, that is:
\begin{longlist}[(2)]
\item[(1)] either of type \textup{{(S)}}, in which case one assumes $\log
f_0 \in C^\alpha$
for some $\alpha> 1$,
\item[(2)] or of type \textup{\hyperref[eqH]{(H)}}, and one assumes $f_0 \in C^\alpha$
for some $1/2<\alpha\le1$.
\end{longlist}
Suppose the prior parameters satisfy (\ref{def-ln}), (\ref{cond-sil})
and (\ref{prk}). Let $\Pi(\cdot\given X^{(n)})$ be the induced
posterior distribution on $\mathcal M_0$.
%Let $\Pi(\cdot| X^{(n)}) \circ\tau_{P_n(j_n)}^{-1}$ be the shifted
%posterior on $\mathcal M_0$ under the transformation $f \mapsto\sqrt
%n (f-\Gamma_n)$.
Then, as $n \to\infty$,
%
%
%e26 #&#
\begin{equation}
\label{iidbvm} \beta_{\mathcal M_0}\bigl(\Pi\bigl(\cdot\given X^{(n)}
\bigr)\circ\tau _{P_n(L_n)}^{-1}, \mathcal N\bigr)
\to^{P^\mathbb N_{0}} 0.
\end{equation}
\end{teo}

%s4 #&#
\section{Some applications}\label{app}

%s4.1 #&#
\subsection{Donsker's theorem for the posterior cumulative distribution function}

Whenever a prior on $f$ satisfies the weak Bernstein--von Mises
phenomenon in the sense of Definition~\ref{weakbvmp}, we can deduce
from the continuous mapping theorem a BvM for integral functionals
$L_g(f) = \int_0^1 g(x) f(x)\,dx$ \textit{simultaneously} for many
$g$'s satisfying bounds on the decay of their wavelet coefficients.
More precisely a bound $\sum_k |\langle g, \psi_{lk} \rangle| \le
c_l$ for all $l$ combined with a weak BvM for $(w_l)$ such that $\sum
c_l w_l<\infty$ is sufficient. Let us illustrate this in a key example
$g_t = 1_{[0,t]}, t \in[0,1]$, where we can derive results paralleling
the classical Donsker theorem for distribution functions and its BvM
version for the Dirichlet process proved in \cite{L83}. With the
applications we have in mind, and to simplify some technicalities, we
restrict to situations where the posterior $f|X$ is supported in $L^2$,
and where the centring $T_n$ in Definition~\ref{weakbvmp} is
contained in $L^2$ (resp., equals $\mathbb{X}^{(n)}$). In this case the
primitives
\begin{eqnarray*}
F(t) &=& \hspace*{-0.5pt}\int_0^tf(x)\,dx,\quad\mathbb
T_n(t)=\hspace*{-0.5pt}\int_0^t
T_n(x)\,dx\quad \biggl(\mbox{resp., } \hspace*{-1pt}\int_0^t
\,dX^{(n)}(x)\biggr),\quad t \in[0,1],
\end{eqnarray*}
define random variables in the separable space $C([0,1])$ of continuous
functions on $[0,1]$, and we can formulate a BvM result in that space.
Different centrings (such as the empirical distribution function) are
discussed below.

%
%th4 #&#
\begin{teo} \label{kolmsmir} Let $\Pi$ be a prior supported in
$L^2([0,1])$, and suppose the weak Bernstein--von Mises phenomenon in
the sense of Definition~\ref{weakbvmp} holds true in $\mathcal M_0(w)$
for some sequence\vspace*{1pt} $(w_l)$ such that $\sum_l w_l 2^{-l/2}<\infty$, and
with centring $T_n$ either equal to $\mathbb X^{(n)}$ or such that
$T_n \in L^2$. For $f \sim\Pi(\cdot|X)$ (conditional on $X$) define
the posterior cumulative distribution function
%
%
%e27 #&#
\begin{equation}
\label{primit} F(t) = \int_0^t f(x)\,dx,\qquad t
\in[0,1].
\end{equation}
Let $G$ be a Brownian motion $(G(t)\dvtx  t \in[0,1])$ in the white noise
model or a $P_0$-Brownian bridge $(G(t) \equiv G_{P_0}(t)\dvtx  t \in
[0,1]), dP_0(x)=f_0(x)\,dx, f_0 \in L^\infty$, in the sampling model.
If $X \sim P_{f_0}$ for some fixed $f_0$, then as $n \to\infty$,
%
%
%e28 #&#
%e29 #&#
\begin{eqnarray}
\label{funlim}
\beta_{C([0,1])}\bigl(\mathcal L\bigl(\sqrt n (F-\mathbb
T_n)\given X\bigr), \mathcal L(G)\bigr)&\to^{P_{f_0}}& 0,
\\
\label{kslim} \beta_\mathbb R \bigl(\mathcal L\bigl(\sqrt n \|F-\mathbb
T_n\|_\infty\given X\bigr), \mathcal L\bigl(\|G\|_\infty\bigr)
\bigr) &\to^{P_{f_0}}&0.
\end{eqnarray}
\end{teo}
\begin{pf}
The mapping
%
%
%e30 #&#
\begin{equation}
\label{defl} L\dvtx \{h_{lk}\} \mapsto L_t\bigl(
\{h_{lk}\}\bigr):=\sum_{l,k}
h_{lk} \int_0^t \psi
_{lk}(x)\,dx,  \qquad t \in[0,1],
\end{equation}
is linear and continuous from $ {\mathcal{M}}_0(w)$ to $L^\infty
([0,1])$ since, for $0<c<C<\infty$,
\begin{eqnarray*}
\biggl\llvert \sum_{l,k} h_{lk} \int
_0^t \psi_{lk}(x)\,dx \biggr\rrvert
&\le&\sum_{l,k} |h_{lk}| \bigl|
\langle1_{[0,t]}, \psi_{lk} \rangle\bigr|
\\
& \le& c\sup_{l,k} \frac{|h_{lk}|}{w_l} \sum
_{l} w_l 2^{-l/2}
\\
&\le& C\| h
\|_{\mathcal M_0},
\end{eqnarray*}
where we have used $\sup_{t \in[0,1]} \sum_k |\langle1_{[0,t]},
\psi_{lk} \rangle| \le c2^{-l/2}$, shown, for example, as in the
proof of Lemma 3 in \cite{GN09}. Also, $L$ coincides with the
primitive map on any function $h \in L^2([0,1])$ with wavelet
coefficients $\{h_{lk}\}\in\ell_2$, since then
\[
L\bigl(\{h_{lk}\}\bigr) = \sum_{l,k}
h_{lk} \langle1_{[0,t]}, \psi_{lk} \rangle= {
\langle}h,1_{[0,t]} {\rangle}= \int_0^t
h(x)\,dx,\qquad t \in[0,1],
\]
in view of Parseval's identity. Moreover, if $\mathbb G$ is a tight
Gaussian random variable in $ {\mathcal{M}}_0$, then the linear
transformation $L(\mathbb G)$ is a tight Gaussian random variable in
$C([0,1])$, equal in law to a Brownian motion or a $P_0$-Brownian
bridge for our choice $\mathbb G=\mathbb W$ or $\mathbb G = \mathbb
G_{P_0}$, respectively, after checking the identity of the
corresponding reproducing kernel Hilbert spaces (cf. \cite{rkhs}, and
using again that $L$ equals the primitive map on $L^2$). The displays
(\ref{funlim})--(\ref{kslim}) now follow from Definition~\ref
{weakbvmp}, the continuous mapping theorem applied to $L$ and $L\circ\|
\cdot\|_\infty^{-1}$, respectively, and noting that $L(f), L(T_n)$
take values in the closed subspace $C([0,1])$ of $L^\infty([0,1])$
under the maintained assumptions. (Although not used here, it can in
fact be checked that the general inclusion $\operatorname{Im}(L)
\subset C([0,1])$ holds true.)
\end{pf}

%
%co1 #&#
\begin{corollary} \label{classd}
Let $\Pi$ be a prior of type \textup{{(S)}} or \textup{\hyperref[eqH]{(H)}}, and suppose the
conditions of Theorem~\ref{thmdens1} are satisfied. Let $F_n(t) =
(1/n)\sum_{i=1}^n 1_{[0,t]}(X_i), t \in[0,1]$, be the empirical
distribution function based on a sample $X_1,\ldots, X_n$ from law
$P_{0}$, and let $F$ be a cumulative distribution function induced by
$\Pi(\cdot|X_1,\ldots, X_n)$ as in (\ref{primit}). Then, as~$n \to
\infty$,
\begin{eqnarray*}
\beta_{L^\infty([0,1])}\bigl(\mathcal L\bigl(\sqrt n (F-F_n)\given X
\bigr), \mathcal L(G_{P_{0}})\bigr) &\to^{P_{0}^\mathbb N}& 0,
\\
\beta_\mathbb R \bigl(\mathcal L\bigl(\sqrt n \|F-F_n
\|_\infty\given X\bigr), \mathcal L\bigl(\|G_{P_{0}}\|_\infty\bigr)
\bigr) &\to^{P^\mathbb N_{0}}& 0.
\end{eqnarray*}
\end{corollary}
\begin{pf}
By Theorems~\ref{thmdens1} and~\ref{kolmsmir} the result is true with
$F_n$ replaced by the primitive $\mathbb T_n$ of $P_n(L_n)$. As in the
proof leading to Remark 9 in \cite{GN09} one shows $\|\mathbb T_n-F_n\|
_\infty= o_P(1/\sqrt n)$, and hence the result follows from the
triangle inequality. (To avoid measurability issues we note that the
result holds for convergence in distribution in $L^\infty([0,1])$ in
the generalised sense of empirical processes (as in \cite{GZ90}), or
in the space of c\`adl\`ag functions on $[0,1]$.)
\end{pf}

Returning to the general setting of Theorem~\ref{kolmsmir}, a natural
credible band for $F$ is to take $C_n, R_n$ such that, with $L$ the map
defined in (\ref{defl}),
%
%
%e31 #&#
\begin{equation}
\label{band0} C_n = \bigl\{F\dvtx  \|F-\mathbb T_n
\|_\infty \le R_n/\sqrt n \bigr\}, \qquad\Pi\circ
L^{-1}(C_n|X)=1-\alpha.
\end{equation}
The proof of the following result implies in particular that $C_n$
asymptotically coincides with the usual Kolmogorov--Smirnov confidence
band. The result is true also with centring $\mathbb T_n = F_n$ (in
which case the proof requires minor modifications related to the
remarks at the end of the proof of Corollary~\ref{classd}).
%
%co2 #&#
\begin{corollary}
Under the conditions of Theorem~\ref{kolmsmir}, let $X \sim P_{f_0},
F_0= \int_0^\cdot f_0(t)\,dt$ and $C_n$ as in (\ref{band0}). Then we
have, as $n \to\infty$,
\[
P_{f_0}(F_0 \in C_n) \to1-\alpha\quad
\mbox{and}\quad R_n \to ^{P_{f_0}} \operatorname{const}.
\]
\end{corollary}
\begin{pf}
The proof is similar to Theorem 1 in \cite{CN13}, replacing $H(\delta
)$ there by $C([0,1])$ (a separable Banach space): the function $\Phi$
in that proof is strictly increasing: any shell $\{g\in C([0,1])\dvtx  s<\|
g\|_\infty<t\}$, $0\le s<t$, contains an element of the RKHS (see
\cite{rkhs}) of Brownian motion [in the case of the white noise model~(\ref{wnmodel})] or of the $P_0$-Brownian bridge (in the case of the
i.i.d. sampling model). Using also Theorem~\ref{CLT} in the sampling
model case, all arguments from the proof of Theorem 1 in \cite{CN13}
go through.
\end{pf}

%re4 #&#
\begin{remark}
Equation (\ref{band0}) [resp., (\ref{def-cs}) below]
reads conditionally on the existence of such a positive real $R_n$.
More generally, one may take a generalised quantile in  \eqref{band0}
[resp., in \eqref{def-cs}]. Then  $C_n$ has credibility $1-\alpha$ asymptotically,
and one can check that the previous corollary [resp., Theorem \ref{confcred}]  continues to hold.
\end{remark}

%s4.2 #&#
\subsection{Confidence bands for $f$}\label{confsec}

Given a posterior distribution $\Pi(\cdot\given X)$ on the parameter
$f$ of a regression or sampling model, we can incorporate the
multiscale approach to construct confidence sets for $f$ in a Bayesian
way. We take an efficient centring $T_n$
[e.g., $\mathbb{X}^{(n)}, P_n(L)$ from above or, when appropriate, the
posterior mean $E(f|X)$] and, given $\alpha>0$ and admissible $w$,
choose $R_n$ and the credible region $C_n$ in such a way that
%
%
%e32 #&#
\begin{equation}
\label{def-cs} C_n = \biggl\{f\dvtx  \sup_{l,k}
\frac{|\langle f- T_n, \psi_{lk}\rangle
|}{w_l} \le\frac{R_n}{\sqrt n} \biggr\}, \qquad  \Pi(C_n|X)=1-
\alpha.
\end{equation}

%
%th5 #&#
\begin{teo} \label{confcred}
Let $w=(w_l)$ be admissible. Suppose the weak Bernstein--von Mises
phenomenon holds true in $\mathcal M_0(w)$ with prior $\Pi$ and
centring $T_n$. Let $C_n$ be as in (\ref{def-cs}). Then, as $n \to
\infty$,
\[
P_{f_0}(f_0 \in C_n) \to1-\alpha,\qquad
R_n \to^{P_{f_0}} \operatorname{const}.
\]
\end{teo}
\begin{pf}
The proof is the same as the one of Theorem 1 in \cite{CN13},
replacing $H(\delta)$ there by $\mathcal M_0(w)$, and using also
Theorem~\ref{CLT} in the sampling model case. %, a separable Banach
%space in which $\ell_2$ is dense.
\end{pf}

The previous theorem can be used to control low frequencies of the
estimation error, and following the multiscale approach one needs to
employ further qualitative information about $f_0$ to control high
frequencies. In the present case, if we assume $f_0 \in C^\gamma$ for
some $\gamma>0$, we can define, for $u_n = w_{j_n}/\sqrt{j_n}$ and
$j_n$ such that $2^{j_n} \sim(n/\log n)^{1/(2\gamma+1)}$, the
confidence set
%
%
%e33 #&#
\begin{equation}
\label{def-cninter} \bar C_n = \bar C_n(\gamma)=C_n
\cap \bigl\{f\dvtx  \|f\|_{C^\gamma} \le u_n \bigr\}.
\end{equation}
The following result combined with Theorems~\ref{thmdens1} and~\ref
{confcred} implies in particular Proposition~\ref{histocov} from the
\hyperref[sec1]{Introduction}.
%
%pr5 #&#
\begin{prop} \label{interpol}
Under the conditions of Theorem~\ref{confcred} suppose $X \sim
P_{f_0}$ where $f_0 \in C^\gamma([0,1])$. Then, with $\bar C_n$ as in
(\ref{def-cninter}), and as $n \to\infty$,
\[
P_{f_0}(f_0 \in\bar C_n) \to1- \alpha\quad
\mbox{and}\quad|\bar C_n|_\infty=O_{P_{f_0}} \bigl( (n/
\log n)^{-\gamma/(2\gamma+1)} u_n \bigr).
\]
\end{prop}
\begin{pf}
For $n$ large enough such that $u_n \ge\|f_0\|_{C^\gamma}$ we have as
$n\to\infty$
\[
P_{f_0}(f_0 \in\bar C_n)
=P_{f_0}(f_0 \in C_n) \to1-\alpha
\]
in view of Theorem~\ref{confcred}. Moreover, for $h=f-g, f,g \in C_n$
arbitrary,
\[
\|h\|_{\mathcal M(w)} \le\|f-T_n\|_{\mathcal M(w)} +
\|g-T_n\| _{\mathcal M(w)} = O \biggl(\frac{R_n}{\sqrt n} \biggr)=
O_{P_{f_0}} \biggl(\frac{1}{\sqrt n} \biggr).
\]
The estimate on $|\bar C_n|_\infty$ now follows from
\[
\|h\|_\infty\le\sum_{l} 2^{l/2}
\max_k \bigl|\langle h, \psi_{lk} \rangle\bigr|
\]
combined with the bound
\begin{eqnarray*}
\sum_{l \le j_n} 2^{l/2}\max_k
\bigl|\langle h, \psi _{lk} \rangle\bigr| &=& \sum_{l \le j_n}
2^{l/2}\sqrt l \frac{w_l}{\sqrt
l} w_l^{-1} \max
_k \bigl|\langle h, \psi_{lk} \rangle\bigr|
\\
&\lesssim&\sqrt{\frac{2^{j_n} j_n}{n}} \frac{w_{j_n}}{\sqrt{j_n}} R_n
\\
&=& O_{P_{f_0}} \biggl( \biggl(\frac{\log n}{n} \biggr)^{\gamma
/(2\gamma+1)}u_n
\biggr),
\nonumber
\end{eqnarray*}
and with
\begin{eqnarray*}
\sum_{l > j_n} 2^{l/2}\max
_k\bigl|\langle h, \psi_{lk} \rangle\bigr| &=& \sum
_{l > j_n} 2^{-l\gamma} 2^{l(\gamma+1/2)}\max
_k\bigl|\langle h, \psi _{lk} \rangle\bigr|
\\
&\le&\|h\|_{C^\gamma} 2^{-j_n\gamma}
\\
&=& O_{P_{f_0}} \biggl( \biggl(
\frac{\log n}{n} \biggr)^{\gamma/(2\gamma+1)}u_n \biggr),
\end{eqnarray*}
completing the proof.
\end{pf}

%
%re5 #&#
\begin{remark}[(Optimal diameter, undersmoothing, adaptation)]\label{opti}
The confidence bands from Propositions~\ref{histocov}
and~\ref{interpol} have diameter equal to the $L^\infty$-minimax rate
over H\"older balls multiplied with an under-smoothing penalty $u_n$,
common in frequentist constructions of confidence bands; see \cite
{H92} and, more recently, \cite{GN10}. If the BvM phenomenon holds for
all admissible sequences $w$ (as in the examples above), then this
sequence can be taken to diverge at an arbitrarily slow rate.

If a quantitative a priori bound $\|f_0\|_{C^\gamma}<B$ is available,
then in the setting of Theorem~\ref{wnbvm} one could use a uniform
wavelet prior [with $\varphi=\1_{[-B,B]}/(2B)$, for some $B>0$]
concentrating on a H\"older ball of radius $B$ (as in Corollary~1,
\cite{CN13}). The set $\bar C_n$ from (\ref{def-cninter}) (even with
$u_n$ replaced by $B$) is then an \textit{exact} level $1-\alpha$
posterior credible set, consisting of the intersection of two
hyper-rectangles in sequence space, and Proposition~\ref{interpol}
applies to give the precise frequentist asymptotics of $\bar C_n$.

% If no bound on $B$ is available we can take $B=B_n =u_n \to\infty$
%without increasing the asymptotic diameter of $\bar C_n$.
%Alternatively, for instance in the case of Gaussian priors, we can use
%posterior information as in Corollary 2 in \cite{CN13} (combined with
%the sup-norm posterior contraction results in Theorem 1 in
%undersmoothed $\bar C_n$.

We can also obtain adaptive confidence bands by using a bandwidth
choice $\hat j_n$ as in \cite{GN10} to estimate $\gamma$ by $\hat
\gamma$ under a self-similarity constraint on $f$, corresponding to an
empirical Bayes-type selection of $\gamma$. More Bayesian approaches
to \textit{adaptive} confidence sets are subject of current research;
see, for example, the recent contribution \cite{SVV13}.
\end{remark}

%s5 #&#
\section{Proofs}\label{proof}

%s5.1 #&#
\subsection{Proof of Theorem \texorpdfstring{\protect\ref{CLT}}{1}}
For $J$ to be chosen below, let $V_J$ be the subspace of $\mathcal
M(w)$ consisting of the scales $l \le J$, and let $\pi_{V_J}(P)$ be
the projection of $f$ onto~${V_J}$. We have by definition of the H\"
older space $C^\gamma$ and assumption
%
%
%e34 #&#
\begin{eqnarray}
\bigl\|P-\pi_{V_{j_n}}(P)\bigr\|_{\mathcal M_0} &=& \sup_{l>j_n}
\frac{\max_k|\langle f, \psi_{lk} \rangle|}{w_l} \nonumber
\\
&\lesssim& w_{j_n}^{-1} 2^{-j_n(\gamma+1/2)}
\\
&=& o\biggl(\frac{1}{\sqrt n} \biggr)\nonumber
\end{eqnarray}
so that this term is negligible in the limit distribution. Writing
$\beta$ for $\beta_{\mathcal M_0}$ and $\sqrt n (P_n(j_n)-\pi
_{V_{j_n}}(P))=\nu_n$, it suffices to show that
%
%
%e35 #&#
\begin{eqnarray}
\beta\bigl(\mathcal L(\nu_n), \mathcal L(\mathbb G_P)
\bigr) &\le&\beta \bigl(\mathcal L(\nu_n), \mathcal L(
\nu_n) \circ\pi_{V_J}^{-1}\bigr)
\nonumber
\\
&&{} + \beta\bigl(\mathcal L(\nu_n) \circ\pi_{V_J}^{-1},
\mathcal L(\mathbb G_P) \circ\pi_{V_J}^{-1}
\bigr)
\\
&&{} + \beta\bigl(\mathcal L(\mathbb G_P), \mathcal L(\mathbb
G_P) \circ \pi_{V_J}^{-1}\bigr)\nonumber
\end{eqnarray}
converges to zero under $P^\mathbb N$. Let $\varepsilon>0$ be given.
The second term is less than $\varepsilon/3$ for every $J$ fixed and
$n$ large enough by the multivariate central limit theorem applied to
\[
\frac{1}{\sqrt n} \sum_{i=1}^n \bigl(
\psi_{lk}(X_i)-E_P\psi _{lk}(X)
\bigr),\qquad k, l \le J,
\]
noting that eventually $j_n >J$. For the first term, by definition of
$\beta$,
%
%
%e36 #&#
\begin{eqnarray}\label{fidap}
&& \beta\bigl(\mathcal L(\nu_n), \mathcal L(
\nu_n) \circ\pi_{V_J}^{-1}\bigr)\nonumber
\\
&&\qquad \le E\bigl\|\sqrt n (
\pi_{V_{j_n}}-\pi_{V_J}) (P_n-P)\bigr\|_{\mathcal M_0}
\\
&&\qquad \le \biggl[ \max_{J < l \le{j_n}} \frac{\sqrt l}{w_l} \biggr] E \max
_{J<l\le{j_n}}l^{-1/2} \max_k \bigl|\bigl
\langle\sqrt n (P_n-P), \psi _{lk}\bigr\rangle\bigr|.\nonumber
\end{eqnarray}
Thus for $J$ large enough this term can be made smaller than
$\varepsilon/3$ if we can show that the expectation is bounded by a
fixed constant. For $M$ a large enough constant, this expectation is
bounded above by $M$ plus
\begin{eqnarray*}
&& \int_M^\infty P \Bigl( \max
_{J \le l \le{j_n}}l^{-1/2} \max_k \bigl|\bigl
\langle\sqrt n (P_n-P), \psi_{lk}\bigr\rangle\bigr| >u \Bigr)\,du
\\
&&\qquad \le\sum_{J \le l \le{j_n}, k} \int_M^\infty
P \bigl( \bigl|\bigl\langle \sqrt n (P_n-P), \psi_{lk}\bigr
\rangle\bigr|>\sqrt{l}u \bigr)\,du
\\
&&\qquad\le\sum_{J \le l \le{j_n}} 2^l \int
_M^\infty e^{-Clu}\,du \lesssim
e^{-C' J M},
\end{eqnarray*}
where\vspace*{1pt} the second inequality follows from an application of Bernstein's
inequality (e.g., \cite{L01}) together with the bounds $P\psi
_{lk}^2\leq\|f\|_{\infty}$ and $\sqrt{l}\|\psi_{lk}\|_\infty\le
\sqrt{l}2^{l/2}=O(\sqrt{n})$ for $l\le j_n$, using the assumption on $j_n$.

For\vspace*{1pt} the third Gaussian term we argue similarly, replacing $\nu_n$ by
$\mathbb G_P$ in (\ref{fidap}) and using that $E \sup_l \max_k
|\mathbb G_P(\psi_{lk})|/\sqrt l<\infty$ by Proposition~\ref{CKR}.

%s5.2 #&#
\subsection{A tightness criterion in $\mathcal M_0$}

The following proposition considers general random posterior measures
$\Pi(\cdot\given X)$ in the setting of Definition~\ref{weakbvmp}.

%
%pr6 #&#
\begin{prop}\label{tightness} Let $\pi_{V_J}, J \in\mathbb N$, be
the projection operator onto the finite-dimensional space spanned by
the $\psi_{lk}$'s with scales up to $l \le J$. Let $f \sim\Pi(\cdot
|X)$, $T_n=T_n(X)$, let $\tilde\Pi_n$ denote the laws of $\sqrt n (f-
T_n)$ conditionally on $X$ and let $\mathcal N$ equal the Gaussian
probability measure on $\mathcal M_0(w)$ given by either $\mathbb W$ or
$\mathbb G_P$ from $P$ with bounded density.

Assume that the finite-dimensional distributions converge, that is,
%
%
%e37 #&#
\begin{equation}
\label{fidij} \beta_{V_J} \bigl(\tilde\Pi_n \circ
\pi_{V_J}^{-1}, \mathcal N \circ\pi_{V_J}^{-1}
\bigr) \to^{P_{f_0}} 0\qquad\mbox{as } n\to \infty,
\end{equation}
and that for some sequence $\bar w=(\bar w_l) \uparrow\infty$, $\bar
w_l/\sqrt l \ge1$,
\begin{eqnarray}\label{inexp}
E \bigl[\|f-T_n\|_{\mathcal M_0(\bar w)}|X\bigr] &=& E \Bigl[\sup
_l \bar w_l^{-1} \max
_{k}\bigl|\langle f - T_n, \psi_{lk}
\rangle\bigr| |X \Bigr]
\nonumber\\[-8pt]\\[-8pt]\nonumber
&=& O_{P_{f_0}} \biggl(\frac{1}{\sqrt n} \biggr).\nonumber
\end{eqnarray}
Then, for any $w$ such that $w_l/\bar w_l \uparrow\infty$ we have, as
$n \to\infty$,
\[
\beta_{\mathcal M_0(w)}(\tilde\Pi_n, \mathcal N) \to^{P_{f_0}} 0.
\]
\end{prop}
%
%
%re6 #&#
\begin{remark}
Inspection of the proof shows that the result still holds true if $f
\sim\Pi(\cdot|X)$ is replaced by $f \sim\bar\Pi(\cdot\given X)$
for random measures $\bar\Pi(\cdot\given X)$ s.t.
\[
\beta_{\mathcal M_0}\bigl(\bar\Pi(\cdot\given X), \Pi(\cdot\given X)\bigr)
\to^{P_{f_0}} 0
\]
as $n \to\infty$. Likewise, the posterior can be replaced by the
conditional posterior $\Pi^{D_n}(\cdot\given X)$ for any sequence of
sets $D_n$ such that
$\Pi(D_n\given X) \to^{P_{f_0}} 1$.
\end{remark}
\begin{pf}
Let us write $\beta= \beta_{\mathcal M_0(w)}$ and decompose
\[
\beta(\tilde\Pi_n, \mathcal N) \le\beta \bigl(\tilde
\Pi_n, \tilde\Pi_n \circ\pi_{V_J}^{-1}
\bigr)+ \beta \bigl(\tilde\Pi_n \circ\pi_{V_J}^{-1},
\mathcal N \circ\pi_{V_J}^{-1} \bigr) + \beta \bigl(\mathcal
N, \mathcal N \circ\pi_{V_J}^{-1} \bigr). %
\]
The second term converges to zero by (\ref{fidij}). The third term
too, arguing as at the end of the proof of Theorem~\ref{CLT} (and
using Proposition~\ref{CKR} or~\ref{CKRP}). For the first term let $f
\sim\Pi(\cdot|X)$ conditional on $X$. Then using (\ref{inexp}) we
can bound the $\beta$-distance by the expectation of the norm and thus by
\begin{eqnarray*}
&& E\bigl[\bigl\|\sqrt n (id-\pi_{V_J}) (f-T_n)
\bigr\|_{\mathcal M(w)}|X\bigr]
\\
&&\qquad \le \biggl[ \sup_{l>J} \frac{\bar{w}_l}{w_l} \biggr] E \Bigl[
\sup_{l>J} \bar{w}_l^{-1} \max
_k \bigl|\sqrt{n}\langle f-T_n, \psi_{lk}
\rangle\bigr| |X \Bigr]
\\
&&\qquad \le\sup_{l>J} \frac{\bar w_l}{w_l} \times
O_{P_{f_0}}(1),
\end{eqnarray*}
which can be made as small as desired for $J$ large enough but fixed.
\end{pf}

%s5.3 #&#
\subsection{Proof of Theorem \texorpdfstring{\protect\ref{wnbvm}}{2}}
We choose integers $j=j_n \to\infty$ such that
%
%
%e38 #&#
\begin{equation}
\label{balance} \sigma_j^{-1}=2^{j(\alpha+1/2)} \sim\sqrt
n\quad\mbox{and note}\quad\sigma_l \lesssim\frac{1}{\sqrt n} \qquad
\forall l > j.
\end{equation}
Conditional on $\mathbb X^{(n)}$, let $f \sim\Pi(\cdot|\mathbb
X^{(n)})$ and, for $\pi_{V_j}$ the projection operator onto~$V_j$,
consider the decomposition in $\mathcal M_0(w)$, under $P_{f_0}$,
\begin{eqnarray*}
\label{zdec} \sqrt n \bigl(f - \mathbb X^{(n)}\bigr) & =& \sqrt n
\bigl(\pi_{V_j}(f) - \pi _{V_j}\bigl(\mathbb X^{(n)}
\bigr)\bigr) + \sqrt n \bigl(f - \pi_{V_j}(f)\bigr)
\\
&&{}+ \sqrt n \bigl(\pi_{V_j}(f_0)-f_0
\bigr) + \bigl(\pi_{V_j}(\mathbb W) - \mathbb W\bigr)
\nonumber
\\
& =& \mathrm{I}+\mathrm{II}+\mathrm{III}+\mathrm{IV}.
\nonumber
\end{eqnarray*}
We verify the conditions\vspace*{1pt} of Proposition~\ref{tightness} above for the
laws $\mathcal L(\sqrt{n}(f-\mathbb{X}^{(n)})|X) = \tilde\Pi_n$
and for the choice $\bar w_l = \sqrt l$. From Theorem 7 in \cite
{CN13}, with Condition 2 verified\vadjust{\goodbreak} in the proof of Theorem 9 of that
paper, we derive condition (\ref{fidij}). Next we verify that (\ref
{inexp}) is satisfied for each of the terms I, II, III, IV, separately.
That is, we check that each term has bounded $ {\mathcal{M}}(\bar
w)$-norm in expectation (and apply Markov's inequality).

(IV) We have as in the proof of Proposition~\ref{CKR} that
\[
E \sup_{k,l} l^{-1/2} \bigl|\mathbb W(
\psi_{lk})\bigr| \le C<\infty.
\]

(III) This term is nonrandom and we have by Condition~\ref{tolk} and
definition of $\sigma_l$, and some constant $0<M<\infty$,
\begin{eqnarray*}
\sqrt n \sup_{l > j,k}l^{-1/2}\bigl\llvert \langle
f_0, \psi_{lk} \rangle \bigr\rrvert & \lesssim& M \sqrt n
\sup_{l > j} l^{-1/2} \sigma_l \lesssim M/
\sqrt j.
\end{eqnarray*}

(II) For $E$ the iterated expectation under $P_{f_0}$ and $\Pi(\cdot
|X)$, we can bound
\[
E \sup_{l > j, k} l^{-1/2} \bigl|\langle f,
\psi_{lk}\rangle\bigr| \le\sum_{l>j}
l^{-1/2} E\max_{k}\bigl|\langle f, \psi_{lk}
\rangle\bigr|.
\]
%
%The expectation of the maxima over $k$ is bounded at the end of the
%proof of Theorem 1 (p.10) in \cite{ic13} ($\delta=1$) by $\sqrt l
%from the one in \cite{ic13}, where an extra downweighting of $
%present choice $\sigma_l= 2^{-l(\al+1/2)}$ as well, with the cut-off
%in \cite{ic13} replaced by the present cut-off $j=j_n$ from (
Denote $f_{lk}:=\langle f, \psi_{lk}\rangle$, $f_{0,lk}:=\langle f_0,
\psi_{lk}\rangle$ and $\varepsilon_{lk}:={\langle}\mathbb{W}, \psi
_{lk} {\rangle}$.
An application of Jensen's inequality yields, for any $t>0$,
\[
E\max_{k}|f_{lk}| \le\frac{1}t \log\sum
_{k} E \bigl(e^{tf_{lk}} + e^{-tf_{lk}}
\bigr).
\]
It is now enough to bound the Laplace transform $E[e^{sf_{lk}}]$ for
$s=t, -t$. Both cases are similar, so we focus on $s=t$,
\begin{eqnarray*}
E\bigl[ e^{t f_{lk}} \bigr] &=& E \frac{ \int e^{t(f_{0,lk}+({v}/{\sqrt{n}}))} e^{-({v^2}/{2})+\varepsilon_{lk}v} ({1}/({\sqrt{n}\sigma_{l}}))
\varphi ( ({f_{0,lk}+({v}/{\sqrt{n}})})/{\sigma_{l}}
 )\,dv }{ \int e^{-({v^2}/{2})+\varepsilon_{lk}v}
({1}/({\sqrt{n}\sigma_{l}}))
\varphi ( ({f_{0,lk}+(v/\sqrt{n})})/{\sigma_{l}}
 ) \,dv }
\\
&=:& E \frac{N_{lk}(t)}{D_{lk}}.
\end{eqnarray*}
To bound the denominator $D_{lk}$ from below, one applies the same
technique as in~\cite{CN13}, proof of Theorem 5. One first restricts
the integral to $(-\sqrt{n}\sigma_{l},\sqrt{n}\sigma_{l})$. Next
one notices, using (\ref{eqP1}), that over this interval the argument of
$\varphi$ lies in a compact set, and hence the function $\varphi$ can
be bounded below by a constant, using (\ref{eqP2}). Next one applies
Jensen's inequality to obtain
\[
D_{lk}\gtrsim e^{-(1/2) \int_{-\sqrt{n}\sigma_{l}}^{\sqrt
{n}\sigma_{l}} (v^2/2) \,dv/(\sqrt{n}\sigma_{l})} \gtrsim e^{-C}.
\]\eject\noindent
To bound the numerator $N_{lk}(t)$ one splits the integral into a part
$N_1$ on $ {\mathcal{A}}:=\{v\dvtx  |f_{0,lk}+v/\sqrt{n}|\le\sigma
_{l}\}$ and a part $N_2$ on its complement $ {\mathcal{A}}^c$. First
\begin{eqnarray*}
E N_1 & \le& e^{t\sigma_{l}} E \int_{ {\mathcal{A}}}
e^{-v^2/2 +
\varepsilon_{lk}v} \varphi \biggl( \frac{f_{0,lk}+ (v/{\sqrt{n}})}{\sigma_{l}} \biggr)\,dv
\\
& \le& e^{t\sigma_{l}} \int_{ {\mathcal{A}}} \varphi \biggl(
\frac
{f_{0,lk}+(v/{\sqrt{n}})}{\sigma_{l}} \biggr)\,dv \leq e^{t\sigma_{l}},
\end{eqnarray*}
using the definition of $ {\mathcal{A}}$ and Fubini's theorem. On
the other hand, the term $N_2$, setting $w=f_{0,lk}+v/\sqrt{n}$ and
using condition~(\ref{eqE}), is bounded by
\begin{eqnarray*}
E N_2 & \le&\int_{(-1,1)^c} e^{t\sigma_{l}w} E\bigl(
e^{-(n/2)(w\sigma_{l}-f_{0,lk})^2+\varepsilon_{lk}\sqrt
{n}(w\sigma_{l}-f_{0,lk})} \bigr) \varphi(w)\,dw
\\
& \le&\int_{(-1,1)^c} e^{t\sigma_{l}w} \varphi(w)\,dw \lesssim
e^{d(\sigma_{l}t)^{2}},
\end{eqnarray*}
for some $d>0$. Conclude, setting $t=\sigma_{l}^{-1}l^{1/2}$, that
\begin{eqnarray*}
E \max_k |f_{lk}| & \leq&\frac{1}t
\log \bigl( 2^l \bigl[ C e^{t\sigma_{l}} + C e^{d(\sigma_{l}t)^{2} } \bigr]
\bigr)
\\
& \lesssim&\frac{l}t + \sigma_{l}+ \frac{1}t (
\sigma_{l}t)^{2} \lesssim \sigma_{l}l^{1/2}.
\end{eqnarray*}
This gives the overall bound,
\[
\sum_{l>j} 2^{-l(1/2+\alpha)} \lesssim2^{-j(1/2+\alpha)}
= O(1/\sqrt n). %
\]

(I) For the frequencies $l \le j_n$ one proves, as in Lemma 1 in \cite
{ic13}, for some constant $C>0$, the sub-Gaussian bound
%
%
%e39 #&#
\begin{equation}
\label{mgf} E_{f_0} E\bigl(e^{t \sqrt n (f_{lk}-X_{lk})}|X\bigr) \le
Ce^{t^2/2}.
\end{equation}
[All that is needed here is $\varphi$ bounded away from zero and
infinity on a compact set, and that $(f_{0,lk}+v/\sqrt n)/\sigma_l$ is
bounded by a fixed constant, true for the $l$'s relevant here.] Then,
by a standard application of Markov's inequality to sub-Gaussian random
variables, writing $\Pr$ for the law with expectation $E_{f_0}E(\cdot
|X)$, we have for all $v>0$ and universal constants $C, C'$ that
\[
\Pr\bigl( \sqrt n |f_{lk}-X_{lk}| >v\bigr) \le C'
e^{-Cv^2}.
\]
We then bound, for $M$ a fixed constant
\begin{eqnarray*}
&& E_{f_0} E \Bigl(\sup_{l \le j} l^{-1/2} \max
_k \sqrt n |f_{lk}-X_{lk}||X \Bigr)
\\
&&\qquad \le M + \int_M^\infty\Pr \Bigl(\sup
_{l \le j,k} l^{-1/2} \max_k \sqrt n
|f_{lk}-X_{lk}|> u \Bigr)\,du. %
\end{eqnarray*}
The tail integral can be further bounded as follows:
\begin{eqnarray*}
&& \sum_{l \le j, k} \int_M^\infty
\Pr \bigl(\sqrt n |f_{lk}-X_{lk}|> \sqrt l u \bigr)\,du
\\
&&\qquad \lesssim\sum_{l \le j}2^l \int
_M^\infty e^{-Cl u^2} \,du \lesssim \sum
_{l \le j} 2^{l} e^{-CM^2l} \le \operatorname{const}
\end{eqnarray*}
for $M$ large enough. This completes the proof of the BvM with
centring \mbox{$T_n = \mathbb{X}^{(n)}$}. From weak convergence toward
$\mathcal N$ of the posterior measures and \mbox{uniform} integrability (as
one can uniformly bound $1+\varepsilon$-moments by the same arguments as
above), we deduce as in Theorem 10 in \cite{CN13} that $\sqrt n
(E(f|X)-\mathbb{X}^{(n)}) \to E\mathcal N=0$ in $\mathcal M_0$ in
probability, so that the posterior mean can replace $\mathbb{X}^{(n)}$
as the centring, completing the proof.

%s5.4 #&#
\subsection{Proof of Theorem \texorpdfstring{\protect\ref{thmdens1}}{3}} \label{prfi}
%A FEW REMARKS ABOUT THIS PROOF, to be fixed:
%
%a) we need to accommodate the new choice of $L_n$ in (\ref{def-ln}) in
%the proofs -- the bias will then be onlY $O(1)$ and not $o(1)$, but
%that is no problem.
%b) I think $j_n$ in the above theorem and the proof below should
%really be $L_n$
%
%---
%
%This needs to be adapted to take the proof from the besov theorem from
%before..
%
%NOTE THAT $\Ga_n$ should now be $P_n(j)$ everywhere.

For $h$ a positive function in $L^2$, denote
\[
c(h) = \log\int_0^1 h(u)\,du,
\]
so that $he^{-c(h)}$ becomes a density on $[0,1]$.
Also, for any element $g$ of $L^2(P_{0})$, denote $\|g\|
_L^2:=P_{0}(g-P_{0}g)^2=\int_0^1 (g-\int_0^1 g)^2 \,dP_{0}$, where $\|
\cdot\|_L$ is a norm on the subspace of $L^2(P_{0})$ consisting of
$P_{0}$-centered functions. For simplicity of notation within the
proof, we denote $X=X^{(n)}$.
%Recall that $\|\cdot\|_\cH$ denotes the $L^2$-norm in the space of
%centered functions with respect to $f_0$.

Let $\rho_n$ the rate in Lemma~\ref{lem-sn}, where we take $M_n=(\log
{n})\wedge
(w_{L_n}/\sqrt{L_n})^{1/2}\to\infty$. For $\varepsilon_n, C$,
respectively, the rate and constant in Lemma~\ref{lem-max}, we set
\[
D_n=\Bigl\{f=e^{T-c(T)}, \|f-f_0
\|_\infty\le\rho_n, \max_{l\le K,
k} \bigl|{\langle}T, \psi_{lk} {\rangle}\bigr|\le C\sqrt{n}\varepsilon_n \Bigr\},
\]
where the part involving the maximum in the definition of $D_n$ is only
needed for the prior \hyperref[eqS2]{(S2)}, and where $K$ is a large enough
integer. Combining Lemmas~\ref{lem-max} and~\ref{lem-sn}, we have
$E_{f_0}\Pi[D_n\given X]\to1$. We also note that for any $l>K$ and
any $k$, the functions $\psi_{lk}$ are orthogonal to constants in $L^2$.
%at $\Ga_n$ of the log-density case. Maybe write lemma to show closeness
%to $P_n$.}

We apply Proposition~\ref{tightness} and the remark after it, with the
posterior conditioned on $D_n$, using the decomposition, for $L=L_n$
and writing $\pi_{V_L}(P_n)$ for $\pi_{V_L}(P_n(L))$,
\begin{eqnarray}
\label{zdec1} \tilde Y_n &=&\sqrt n \bigl( f - P_n(L)
\bigr)
\nonumber
\\
& =& \sqrt n \bigl( \pi_{V_L}(f) - \pi_{V_L}(P_n)
\bigr) + \sqrt n \bigl(f - \pi _{V_L}(f)\bigr)
\nonumber
%+ \sqrt n (\pi_{V_L}(P_n)-P_n(L)) \nonumber\\
=:
Y_n + r_n.
\end{eqnarray}
%
%where $Y_n$ has law $\bar\Pi(\cdot|X)$.
Thus to prove (\ref{iidbvm}) it suffices to show (i) that $\tilde Y_n
- Y_n$ is asymptotically negligible and to check the conditions of
Proposition~\ref{tightness}, that is, (ii) that
(\ref{inexp}) holds for $Y_n$, and (iii) that finite-dimensional
convergence (\ref{fidij}) occurs. %For our choice of $L$ we have
%(c)$=0$.
\begin{longlist}[(iii)]
\item[(i)] The term $r_n$ is zero in the case of the histogram prior \hyperref[eqH]{(H)},
by definition of the prior and orthogonality of the Haar basis.
To check that $r_n$ is negligible for the log-density priors {(S)},
let us write $f=f_0+(f-f_0)$ and study separately $\pi_{V_L^c}f_0$ and
$\pi_{V_L^c}(f-f_0)$. %TO BE adapted !!! Let $L=j_n$ be the largest
%integer such that $2^j\le n^{1/(1+2\al)}$.
For both choices of $L$, we have $L\ge l_n$, so
\begin{eqnarray*}
\sqrt{n}\|\pi_{V_L^c}f_0\|_{\mathcal M_0} & \le&\sqrt{n}
\sup_{l>l_n} w_l^{-1} \max
_{k} \bigl| {\langle}f_0,\psi_{lk} {
\rangle}\bigr|
\\
& \lesssim&\sqrt{n}\bigl(l_n^{1/2}/w_{l_n}\bigr)
\sup_{l>l_n} l^{-1/2} 2^{-l((1/2)+\alpha)} = o(1),
\end{eqnarray*}
using that $f_0\in C^{\alpha}$, admissibility of $w$ and the
definition of $l_n$. Also,
\begin{eqnarray*}
&& \sqrt{n}\int\bigl\|\pi_{V_L^c}(f-f_0)\bigr\|_{\mathcal M_0}\,d\Pi
^{D_n}(f\given X)
\\
&&\qquad = \sqrt{n} \int\sup_{l>L} w_l^{-1}
\max_{k} \bigl| {\langle}f-f_0,\psi_{lk}
{\rangle}\bigr| \,d\Pi^{D_n}(f\given X)
\\
&&\qquad \le\sqrt{n} \sup_{l>L} {w_l}^{-1} \|
\psi_{lk}\|_1 \int\|f-f_0 \|_\infty
\,d\Pi^{D_n}(f\given X)
\\
&&\qquad \lesssim\sqrt n \bigl(L^{1/2}/w_L\bigr) L^{-1/2}
2^{-L/2} M_n \bigl(2^L L/n\bigr)^{1/2}
=o(1),
\end{eqnarray*}
using $\|\psi_{lk}\|_1 \lesssim2^{-l/2}$ and Lemma~\ref{lem-sn} with
$M_n\to\infty$
as defined above.

\item[(ii)] To control $Y_n$, a key ingredient is a bound on the following
exponential moment restricted to $D_n$. Below we prove that for
universal constants $c_1, c_2$ and $|s|\le\sqrt{l}$, for any $l\le L$
and $k$,
%
%
%e40 #&#
\begin{equation}
\label{eq-bound} \int e^{ s \sqrt{n}{\langle}f-P_n, \psi_{lk} {\rangle}} \,d\Pi ^{D_n}(f\given X) \le
c_1 e^{c_2 s^2}\Pi(D_n\given X)^{-1}.
\end{equation}
%
%A consequence of \eqref{eq-bound} is that every coordinate of the
%projection
%$\pi_{V_j}f$ concentrates at rate $\rn$ under the posterior (on $D_n$)
%around its efficient estimator given by \eqref{gan-sob}.
Suppose for now that (\ref{eq-bound}) is established. Then aiming at
checking (\ref{inexp}) with \mbox{$\bar w_l=\sqrt{l}$}, we can use it in the
study of%, denoting $\zeta:=\zeta_{n,l,k}(f)=\rn\psg f-P_n,\psi_{lk}
\[
\sqrt{n}\bigl\| \pi_{V_L}(f-P_n) \bigr\|_{ {\mathcal{M}}_0(\sqrt{l})}= \sqrt{n}
\max_{l\le L} l^{-1/2}\max_k \bigl|{
\langle}f-P_n,\psi_{lk}{\rangle}\bigr|
\]
in expectation under $\Pi^{D_n}(\cdot\given X)$.
Denoting $E$ and $\operatorname{Pr}$, respectively, for expectation and
probability under $\Pi^{D_n}(\cdot\given X)$, for any $M>0$, with
$ {\mathcal{M}}_0= {\mathcal{M}}_0(\sqrt{l})$,
\begin{eqnarray*}
\sqrt{n}E \bigl\| \pi_{V_L}(f-P_n) \bigr\|_{ {\mathcal{M}}_0} & \le& M +
\int_M^\infty\operatorname{Pr} \bigl[\sqrt{n}\bigl\|
\pi_{V_L}(f-P_n) \bigr\| _{ {\mathcal{M}}_0} >u \bigr]\,du
\\
& \le& M+ \sum_{l<L, k} \int_M^\infty
\operatorname{Pr} \bigl[ \sqrt {n}\bigl|{\langle}f-P_n,
\psi_{lk}{\rangle}\bigr|>\sqrt{l} u \bigr]\,du.
\end{eqnarray*}
An application of Markov's inequality for $u>0$ leads to
\begin{eqnarray*}
\operatorname{Pr} \bigl[ \sqrt{n}\bigl|{\langle}f-P_n,
\psi_{lk}{\rangle }\bigr|> \sqrt{l} u \bigr] & \le& e^{-l u} E \bigl[
e^{ \sqrt{l n} |{\langle}f-P_n,\psi
_{lk}{\rangle}|} \bigr].
\end{eqnarray*}
Combining the last two bounds with (\ref{eq-bound}) leads to
\begin{eqnarray*}
\sqrt{n}E \bigl\| \pi_{V_L}(f-P_n) \bigr\|_{ {\mathcal{M}}_0} &
\lesssim &M + \sum_{l<L} 2^l
e^{c_2 l} \Pi(D_n\given X)^{-1} \int
_M^\infty e^{- l u}\,du
\\
& \lesssim &M+\Pi(D_n\given X)^{-1}\sum
_{l<L} l^{-1} e^{ [l\log{2}+ lc_2 -l M ]}.
\end{eqnarray*}
For $M$ large enough, the last display is bounded by $M+C\Pi(D_n\given
X)^{-1}$. Since $\Pi(D_n\given X)^{-1}\to1$ in probability, one obtains
%
%
%e41 #&#
\begin{equation}
\label{interm} \sqrt{n}E \bigl\| \pi_{V_L}(f-P_n)
\bigr\|_{ {\mathcal{M}}_0(\sqrt{l})} = O_{P_{f_0}}(1).
\end{equation}
Combining (\ref{interm}) with Markov's inequality and Proposition~\ref
{tightness}, we see that the BvM result will follow from (\ref
{eq-bound}), and from convergence of finite-dimensional distributions
that we check in point (iii) below.

Now we check (\ref{eq-bound}), in two steps. First, Lemma~\ref
{lem-lapu} below enables us to incorporate the term ${\langle
}f-P_n,\psi_{lk} {\rangle}$ into the likelihood coming from Bayes'
formula applied to $d\Pi^{D_n}(f\given X)$ and reduces the problem to a
change of measure with respect to the prior. Next, this change of
measure is handled below.

Let us now apply Lemma~\ref{lem-lapu} below to $\Pi_n=\Pi^{D_n}$ for
$\Pi$, one of the considered priors. Set $\gamma_n=\psi_{lk}$. First
note that $\|\tilde\gamma_n\|_L^2 = \int\psi_{lk}^2 f_0\lesssim\|
f_0\|_{\infty}$, which is bounded by assumption for $f_0\in
{\mathcal{F}}_0$. Next note that, for $l\le L$,
\[
\|\tilde\gamma_n\|_\infty\lesssim\|\psi_{lk}
\|_\infty\lesssim 2^{l/2} \lesssim2^{L/2}.
\]
For $h$ the Hellinger distance we have $h(f,f_0)^2 \lesssim\|f-f_0\|
_2^2\le\|f-f_0\|_\infty^2$ valid for $f_0\in{\mathcal{F}}_0$.
Hence\vspace*{1pt} on $D_n$ we have that $h(f,f_0)\lesssim\rho_n$. Since $\alpha
>1/2$, we have $2^{L/2}\log{n} \le\rho_n^{-1}$, so one can apply
Lemma~\ref{lem-lapu} with $a_n=\rho_n$ and deduce %writing $\Pi_n$ as
%a shorthand for $\Pi[\cdot\given X^{(n)}]$,
%
%e42 #&#
\begin{equation}
\label{eq-inter} \int e^{ s \sqrt{n}{\langle}f- P_n, \psi_{lk} {\rangle}} \,d\Pi ^{D_n}(f\given X) \lesssim
\frac{e^{Cs^2}}{\Pi(D_n\given X)} \frac{ \int_{D_n}
e^{\ell_n(f_s)-\ell_n(f_0) } \,d\Pi(f) }{
 \int e^{\ell_n(f)-\ell_n(f_0)}\,d\Pi(f) }.
\end{equation}

Now we are ready to change variables in the last ratio. For each of the
examples of priors considered, we show that this ratio is bounded from
above by a constant as $n\to\infty$.

We start with case {(S)}. By definition, the quantity $f_s$ is a
function of $\log f - s\tilde\gamma_n/\sqrt{n}$. Next, notice that
any constant in this expression vanishes due to the subtraction of the
renormalising constant $c(\cdot)$. In particular, the expression is a
\mbox{function} of $T - s \gamma_n/\sqrt{n}$, where $T$ is defined in (\ref
{def-T}). The law of $T$ is induced by a finite product of probability
measures, via the distributions of the coordinates of $T$ over $\{\psi
_{lk}\}$ with $l<L$.
% (in the case {\bf(GCb)} one sees the prior as a countable product
%measure over the coordinates of the eigenvector basis $\{e_l\}$).
Since $\gamma_n=\psi_{lk}$, only one coordinate of the product
measure defining $T$ is affected by the subtraction of $s\gamma
_n/\sqrt{n}$. The next step is to change variables in the numerator of
the ratio above by shifting the corresponding coordinate by $s/\sqrt{n}$.

For \hyperref[eqS1]{(S1)}, the change in density on this coordinate
can be measured by
%
%
%e43 #&#
\begin{equation}
\label{eq-rat} \frac{\varphi_H({\cdot}/{\sigma_{l}})}{\varphi_H(({\cdot-s/\sqrt{n}})/{\sigma_{l}})},
\end{equation}
whose logarithm is bounded above in absolute value by $1/(\sqrt
{n}\sigma_{l}) \lesssim1$, since by assumption, $\log\varphi_H$ is
Lipschitz and using (\ref{cond-sil}) combined with $l \le L$.

In case \hyperref[eqS2]{(S2)}, the prior on each coordinate is Gaussian, and
if $\theta_{lk}$ denotes the integrating variable with respect to the
coordinate $l,k$ (corresponding to integrating out the law of ${\langle
}T,\psi_{lk}{\rangle}$) in the considered ratio of integrals, we have
%
%
%e44 #&#
\begin{eqnarray}
\log\frac{\varphi_G({\theta_{lk}}/{\sigma_{l}})}{\varphi
_G((\theta_{lk}-s/\sqrt{n})/{\sigma_{l}})} & =& \frac{1}{n\sigma_{l}^2} - \frac{s}{\sqrt{n}\sigma_{l}^2}\theta
_{lk}. \label{cv-gauss} %\\
%& = \frac{1}{n\sil^2} - \frac{s}{\rn\sil^2}(\te_{lk}-\te_{0,lk}) -
\end{eqnarray}
Recall that we work on the set $D_n$, on which we have the following
inequalities: $\|\log(f/f_0)\|_2\leq\|\log(f/f_0)\|_\infty\lesssim
\rho_n$, using that $f_0$ is bounded from below.\vadjust{\goodbreak} Moreover, note that
by definition of $T$, and if $g:=\log f$, $g_0:=\log f_0$, it holds
${\langle}g-g_0, \psi_{lk} {\rangle}={\langle}T-c(T)-g_0, \psi
_{lk}{\rangle}$. Since $c(T)$ is a constant, and $\psi_{lk}$ are
orthogonal to constants for $l\ge K$, $K$ large enough,\vspace*{1.5pt} we deduce, if
$g_{0,lk}:=
{\langle}g_0,\psi_{lk}{\rangle}$, that on $D_n$ we have $(\theta
_{lk}-g_{0,lk})^2\lesssim\rho_n^2$, as soon as $l\ge K$. So, for
$K\le l\le L$, we have
\begin{eqnarray*}
\bigl|(\ref{cv-gauss})\bigr| & \lesssim&\frac{1}{n\sigma_{l}^2} + \frac{\rho
_n|s|}{\sqrt{n}\sigma_{l}^2} +
\frac{|g_{0,lk}||s|}{\sqrt{n}\sigma
_{l}^2} \lesssim1+ \frac{\rho_n}{\sqrt{n}\sigma_{l}^2}\sqrt{l} + \frac{|g_{0,lk}|}{\sigma_{l}}
\sqrt{l}.
\end{eqnarray*}
Since $\log f_0$ belongs to $C^\alpha$ by assumption and with (\ref
{cond-sil}), the last term in the last display is at most a constant.
We also have $\sqrt{l}\rho_n\lesssim\sqrt{n}\sigma_{l}^2$ using
(\ref{cond-sil}) in the Gaussian case, thus the previous display is at
most a constant on $D_n$. Now we are left with the indexes such that
$l\le K$. For those,\vspace*{1pt} by definition of the set $D_n$, (\ref{cv-gauss})
is in absolute value less than $(n^{-1}+\varepsilon_n)\sqrt{l}\sigma
_{l}^{-2}$. Since $l\le K$ with $K$ fixed, the last expression is
bounded, which yields (\ref{eq-bound}).

Finally, the case of the histogram prior \hyperref[eqH]{(H)} is treated by
studying the effect of the change of variables on the Dirichlet
distribution. The argument is similar to~\cite{ic13}, Section~4.4 and
is omitted.

\item[(iii)] Convergence of finite-dimensional distributions (\ref{fidij}).
This can be seen to consist of establishing BvM results for the
projected law of the posterior distribution on any fixed
finite-dimensional subspace $V=\operatorname{Vect}\{ \psi_{lk},
(l,k)\in{\mathcal{T}}\}$, with $ {\mathcal{T}}$ a \emph{finite} admissible set of indexes. By Cram\'er--Wold, this is the same
as showing a BvM for estimating\vspace*{1pt} the linear functional ${\langle}f,\psi
_{ {\mathcal{T}}}{\rangle}_2$, with $\psi_{ {\mathcal
{T}}}:=\sum_{(l,k)\in{\mathcal{T}}} t_{l,k}\psi_{lk}$ and
$t_{l,k}\in\mathbb{R}$. Denote by $\pi_ {\mathcal{T}}$
%:\cB\to\RR$
the mapping, %the orthogonal projection onto $V$,
for any finite set of indices $ {\mathcal{T}}$,
\[
\pi_ {\mathcal{T}}\dvtx f\to{\langle}f,\psi_ {\mathcal
{T}}{\rangle}_2.
\]
Then it is enough to show that, for any finite $ {\mathcal{T}}$,
\[
\beta_{V} \bigl(\Pi(\cdot\given X)\circ\tau^{-1}_{P_n}
\circ\pi _ {\mathcal{T}}^{-1}, N\bigl(0,\|\psi_ {\mathcal{T}}
\|_L^2\bigr) \bigr)\to0,
\]
as $n\to\infty$. Since $ {\mathcal{T}}$ is finite, the
supremum-norm $\|\psi_ {\mathcal{T}}\|_{\infty}$ is bounded. Thus
the techniques of \cite{cr13} can be used for the considered priors.

%In the case of histogram priors \hyperref[eqH]{(H)}, the previous display follows
%from Theorem~4.2 in \cite{cr13} applied to dyadic histograms: one
%considers the linear \mbox{functional} $f\to{\langle}f,\psi_{ {\mathcal
%{T}}}{\rangle}_2$ so $\tilde r=0$ in that theorem, and it suffices to
%check conditions \mbox{(4.9)--(4.10)} with $K_n=L_n$. Condition (4.9) boils
%down, with $g_{[K_n]}$ the $L^2$ projection of a given function $g$ in
%$L^2$ onto the space of regular dyadic histograms of level~$K_n$, to
%$\sqrt{n}\int(\psi_{ {\mathcal{T}}}-\psi_{ {\mathcal
%{T}},[K_n]})(f_{0,[K_n]}-f_0)=o(1)$. But $\psi_{ {\mathcal{T}}}$
%is a dyadic histogram of fixed meshwidth, thus $\psi_{ {\mathcal
%{T}},[K_n]}= \psi_{ {\mathcal{T}}}$ for large enough $n$, since
%$K_n=L_n\to\infty$, so (4.9) trivially holds. Finally (4.10) in \cite
%{cr13} follows from the fact that $K_n=L_n\to\infty$.

In the case of histogram priors (H), the previous display follows from the section on random histograms in \cite{cr13}, applied to dyadic histograms. The functional $\pi_\cT$ above is linear, so the no-bias condition in \cite{cr13} amounts to check, with $g_{[K_n]}$ the $L^2$ projection of a given function $g$ in $L^2$ onto the space of regular dyadic histograms of level $K_n$, that $\rn \int (\psi_{\cT}-\psi_{\cT,[K_n]})(f_0-f_{0,[K_n]})=o(1)$.  But $\psi_{\cT}$ is a dyadic histogram of fixed meshwidth, thus $\psi_{\cT,[K_n]}= \psi_{\cT}$ for large enough~$n$, since $K_n=L_n\to\infty$, so this trivially holds. Finally,\vspace*{1pt} since $K_n=L_n\to\infty$, the variances $\int  (\psi_{\cT,[K_n]}-\int \psi_{\cT,[K_n]}f_0)^2 f_0$ converge to $\int  (\psi_{\cT}-\int \psi_{\cT}f_0)^2 f_0$.

In the case of log-density priors {(S)}, one applies the general result on density estimation in \cite{cr13} (Theorem 4.1). The set $A_n$ in that statement should be replaced by the set $D_n$ defined above. Since $D_n$ is contained in $A_n=\{f\dvtx\ \|f-f_0\|_1\le \rho_n\}$ and $\Pi(D_n\given X)$ tends to $1$ in probability, the proof of that Theorem goes through without further changes. It thus suffices to verify  that the ratio of integrals in the former theorem from  \cite{cr13} holds when the functional $f\to \psg f,\psi_{\cT}\psd_2$ is considered.\vadjust{\goodbreak}
Note that this is the same as proving that the ratio on the right-hand side of \eqref{eq-inter} goes to $1$, with $\psi_{lk}$ replaced by $\psi_{\cT}$ and now $f_s=f e^{-t\psi_{\cT}-c(f e^{-t\psi_{\cT}})}$. Since only a {\em finite} number
 of $\psi_{lk}$s are involved in the sum defining $\psi_{\cT}$, the ratio involved in the change of variables tends to $1$ in probability: in the case of log-Lipschitz densities, one uses a finite number of times the  bound $1/(\rn\sil)$ for the logarithm of \eqref{eq-rat}, which is of the order $1/\rn$ because $l$ is now bounded. For the Gaussian density case, one argues similarly.
\end{longlist}

We conclude with the following auxiliary results: for $L_n=l_n$ these
are Lemmas~3, 6, 9 and Theorems 2, 3 in \cite{ic13}, and the case $L_n
= j_n$ is proved in the same way. Let $h(f,g)$ denote the Hellinger
distance between two given densities $f,g$, and write $\ell_n(f) =
(1/n) \sum_{i=1}^n \log f(X_i)$ for $f>0$.

%
%le1 #&#
\begin{lem} \label{lem-lapu}
Let $f_0$ belong to $ {\mathcal{F}}_0$. %=\cF_{m,M}$.
Let $\{a_n\}$ be a sequence of reals such that $n a_n^2\ge1$ for any
$n\ge1$.
Let $\{\Pi_n\}$ be a collection of priors on densities restricted to
the set
$\{f, h(f,f_0)\le a_n\}$.
Let $\{\gamma_n\}$ be an arbitrary sequence in $L^{\infty}[0,1]$. Set
$\tilde\gamma_n:= \gamma_n - P_{0}\gamma_n$. Suppose, for some
$m>0$ and %some $C_2>0$ and
all $n\ge1$,
\[
\| \tilde\gamma_n \|_L \le m, \qquad\| \tilde
\gamma_n \|_{\infty
} \le\bigl(4a_n\log(n+1)
\bigr)^{-1}.
\]
Then there exist $C>0$ depending on $m, \|f_0\|_{\infty}$ only %$N_0
such that for any $n\ge1$ and $|t|\le\log{n}$, with $W_n(\gamma
_n)=\sqrt{n}(P_n-P_0)\gamma_n$,
\[
E^{\Pi_n} \bigl[ e^{ t \sqrt{n}{\langle}f-f_0, \gamma_n {\rangle
}_2 } \given X^{(n)} \bigr] \le
e^{ Ct^2 + tW_n(\gamma_n)} \frac{ \int e^{\ell_n(f_t)-\ell_n(f_0) } \,d\Pi_n(f)
}{\int e^{\ell_n(f)-\ell_n(f_0)}\,d\Pi_n(f) },
\]
where $f_t$ is defined by $\log f_t = \log f - t\tilde\gamma_n/\sqrt
{n}- c(f e^{- t\tilde\gamma_n/\sqrt{n}})$.
\end{lem}

%
%le2 #&#
\begin{lem} \label{lem-tech}
Let $f,f_0$ be two densities such that $f_0$ is bounded away from
infinity. % and
%$\|f-f_0\|_1\le1/2$ $h(f,f_0) \le1/2$.
Let $g$ be an element of $L^{\infty}$ such that $h(f,f_0)\|g\|_{\infty
} \le C_1$ and $\|g\|_2 \le C_2$, for some constants $C_1, C_2>0$. Then
\[
\bigl|(P-P_0)g^2\bigr| %\left| \int_0^1 g_l^2(f-f_0) \right|
\le C_1^2+C_1
\sqrt{4C_2\|f_0\|_{\infty}+
C_1^2}.
\]
\end{lem}

%
%le3 #&#
\begin{lem} \label{lem-max}
Let $f_0\in{\mathcal{F}}_0$, and suppose $\log f_0 \in C^\alpha
$, $\alpha>1$. Let $\varphi=\varphi_G$ and $\sigma_{l}$ satisfy
(\ref{cond-sil}), $L_n$ be as defined in (\ref{def-ln}) and the prior
$\Pi$ be defined by (\ref{def-prior-log}). Then there exists $\nu>0$
such that if $\varepsilon_n=(\log n)^\nu n^{-\alpha/(2\alpha+1)}$,
for $C>0$ large enough and
any fixed given integer $K$,
\[
E_{f_0}\Pi \Bigl[ \max_{\lambda\le K,\mu} \bigl|{\langle}T, \psi
_{\lambda\mu} {\rangle}\bigr| \le C\sqrt{n}\varepsilon_n\given
X^{(n)} \Bigr] \to1.
\]
\end{lem}

Finally, for any $\alpha>0$ and $n\ge2$, let us set
$\varepsilon_{n,\alpha}^*:= (n/\log n)^{-\alpha/(2\alpha+1) }$.

%
%le4 #&#
\begin{lem} \label{lem-sn}
Let $\Pi$ be of the type \textup{{(S)}} or \textup{\hyperref[eqH]{(H)}} with $L_n$ as in (\ref{def-ln}) and suppose (\ref{cond-sil}) and (\ref{prk}) are,
respectively, satisfied for the corresponding prior.
Suppose $f_0$ belongs to $C^\alpha$, $1/2<\alpha\le1$ in the case of
prior \textup{\hyperref[eqH]{(H)}} and $\log f_0\in C^\alpha$, $\alpha>1$ for\vadjust{\goodbreak} priors \textup{{(S)}}.
Then, as $n\to\infty$,
%
%
%e45 #&#
\begin{equation}
\label{rate-sn} E_{f_0}\Pi\bigl[ f\dvtx  \|f-f_0
\|_\infty> \rho_n \given X^{(n)} \bigr] \to0,
\end{equation}
where, for an arbitrary sequence $M_n\to\infty$, $\rho_n^2=M_n^2 L_n
2^{L_n}/n$. That is, $\rho_n= M_n\varepsilon_{n,\alpha}^*$ if
$L_n=l_n$ and
$\rho_n=M_n\varepsilon_{n,\alpha}^* \sqrt{\log n}$ if $L_n=j_n$.
\end{lem}

\section*{Acknowledgments}
The authors thank two anonymous referees and
the Associate Editor for comments that helped to improve the presentation of the
paper. Richard Nickl is grateful to the LPMA at Universit\'e Paris VII Denis
Diderot for its hospitality during a visit in October--November~2012 where this
research was initiated.

%suskaldyti doi

% imsref loaded by linak, 2014-07-03 12:14:34
%
% imsref loaded by linak, 2014-07-04 11:45:05

\printaddresses
\end{document}